\documentclass{article}
\usepackage{amsfonts,amssymb}
\usepackage{amsthm}
\usepackage[centertags]{amsmath}
\def\bf{\bfseries}

\newcommand{\eps}{\varepsilon}

\newtheorem{theorem}{Theorem}
\newtheorem{itlemma}{Lemma}[section]
\newtheorem{itclaim}{Claim}[section]

\newtheorem{itproposition}[itlemma]{Proposition}
\newtheorem{itcorollary}[itlemma]{Corollary}
\newtheorem{itremark}[itlemma]{Remark}
\newtheorem{itdefinition}[itlemma]{Definition}
\newtheorem{itexample}[itlemma]{Example}
\newenvironment{claim}{\begin{itclaim}\rm}{\end{itclaim}}
\newenvironment{lemma}{\begin{itlemma}\rm}{\end{itlemma}}
\newenvironment{remark}{\begin{itremark}\rm}{\end{itremark}}
\newenvironment{corollary}{\begin{itcorollary}\rm}{\end{itcorollary}}
\newenvironment{proposition}{\begin{itproposition}\rm}{\end{itproposition}}
\newenvironment{definition}{\begin{itdefinition}\rm}{\end{itdefinition}}
\newenvironment{example}{\begin{itexample}\rm}{\end{itexample}}

\newcommand{\bcl}[1]{\begin{claim}\label{#1}}
\newcommand{\bl}[1]{\begin{lemma}\label{#1}}
\newcommand{\br}[1]{\begin{remark}\label{#1}}
\newcommand{\bt}[1]{\begin{theorem}\label{#1}}
\newcommand{\bd}[1]{\begin{definition}\label{#1}}
\newcommand{\bp}[1]{\begin{proposition}\label{#1}}
\newcommand{\bc}[1]{\begin{corollary}\label{#1}}
\newcommand{\bfact}[1]{\begin{fact}\label{#1}}
\newcommand{\bex}[1]{\begin{example}\label{#1}}
\newcommand{\bem}[1]{\begin{example}\label{#1}}
\newcommand{\ec}{\end{corollary}}
\newcommand{\eex}{\end{example}}
\newcommand{\eem}{\end{example}}
\newcommand{\el}{\end{lemma}}

\newcommand{\er}{\end{remark}}
\newcommand{\et}{\end{theorem}}
\newcommand{\ed}{\end{definition}}
\newcommand{\ep}{\end{proposition}}
\newcommand{\epr}{\end{proof}}
\newcommand{\bpr}{\begin{proof}}
\newcommand{\ecl}{\end{claim}}

\newcommand{\beq}{\begin{eqnarray}}
\newcommand{\eeq}{\end{eqnarray}}
\newcommand{\beqn}{\begin{eqnarray*}}
\newcommand{\eeqn}{\end{eqnarray*}}
\newcommand{\bi}{\begin{itemize}}
\newcommand{\ei}{\end{itemize}}
\newcommand{\ben}{\begin{enumerate}}
\newcommand{\een}{\end{enumerate}}

\newcommand{\R}{{\mathbb R}}  
\newcommand{\N}{{\mathbb N}}  

\newcommand{\url}[1]{\mathtt{#1}}

\def\Id{{1\kern-.4em 1}}

\def\epi{{\text{epi}\;}}

\def\C{{\cal C}}
\def\cv{{\cal V}}
\title{Infinitesimal Characterizations for Strong Invariance and Monotonicity
for Non-Lipschitz Control Systems\footnote{The author thanks
Eduardo Sontag for suggesting the problems addressed in this work
and Peter Wolenski for commenting on an earlier draft.
 Supported
 by  Louisiana Board of Regents Support Fund Grant
 LEQSF(2003-06)-RD-A-12.}}
\author{
Michael Malisoff\\ Department of Mathematics\\ Louisiana State
University\\ Baton Rouge, LA 70803-4918\\ \tt{malisoff@lsu.edu}
}\date{
\today}
\begin{document}
\maketitle

\begin{abstract}
We provide new infinitesimal characterizations for strong
invariance of multifunctions  in terms of Hamiltonian inequalities
and tangent cones.  In lieu of the standard local Lipschitzness
assumption on the multifunction, we assume a new feedback
realizability condition that can in particular be satisfied by
control systems that are discontinuous in the state variable. Our
realization condition is based on H. Sussmann's unique limiting
property, and allows a more general class of feedback realizations
than is allowed by the recent strong invariance characterizations
\cite{KMW04}.
 We also
give new nonsmooth monotonicity characterizations for control
systems that may be discontinuous in the state.\medskip

\end{abstract}

\section{Introduction}
\label{sec1} The theory of flow invariance plays an important role
in much of modern control theory and optimization (see \cite{AS03,
CLSW98, CS03, FPR95, V00, WZ98}). For a given set valued dynamics
$F$ evolving on $\R^n$ and a subset $S\subseteq\R^n$, the theory
provides necessary and sufficient conditions under which  the pair
$(F,S)$ is {\em strongly invariant}, meaning,  for each $T>0$ and
each trajectory $y:[0,T]\to\R^n$ of $F$ starting at a point in $S$
we have $y(t)\in S$ for all $t\in [0,T]$.
   For the special
case where $F$ is locally Lipschitz and $S\subseteq\R^n$ is
closed, infinitesimal characterizations for strong invariance are
well known.  For example, if $F:\R^n\rightrightarrows \R^n$ is
locally Lipschitz and nonempty, compact, and convex valued with
linear growth and $S\subseteq\R^n$ is closed,  then it is well
known (cf. \cite[Chapter 4]{CLSW98}) that $(F,S)$ is strongly
invariant if and only if $F(x)\subseteq T^C_S(x)$ for all $x\in
S$, where $T^C_S$ denotes the Clarke tangent cone (cf. Section
\ref{sec2} and Appendix \ref{nonsmo} for the relevant
definitions). However, this cone characterization can fail if $F$
is non-Lipschitz as illustrated in the following simple example:
Take $n=1$, $S=\{0\}$, $F(0)=[-1,+1]$, and $F(x)=\{-{\rm
sign}(x)\}$ for $x\ne 0$.
 Then $T^C_S(0)=\{0\}$, so $F(0)\not\subseteq T^C_S(0)$.  However,
 $(F,S)$ is strongly invariant.  This example  is covered by the
main sufficient conditions for strong invariance in \cite{KMW04,
RW03}.

On the other hand, consider the controlled differential inclusion
\begin{equation}
\label{mg} \dot x\in G(x,\alpha):=\prod_{i=1}^ng_i(x,\alpha_i)\,
D_i(x)
\end{equation}
where each factor $g_i:\R^n\times A\to \R: (x,a)\mapsto g_i(x,a)$
is locally Lipschitz, $A\subseteq\R^m$ is compact, $\alpha_i\in
{\cal A}:=\{{\rm measurable\ } [0,\infty)\to A\}$ for each $i$,
$\alpha=(\alpha_1,\alpha_2,\ldots, \alpha_n)$, and each of the
multifunctions $D_i:\R^n\rightrightarrows \R$ is Borel measurable.
(Throughout this note, $\prod_{i=1}^n S_i=S_1\times \cdots \times
S_n$ for subsets $S_i\subseteq \R$ and $\prod_{i=1}^n
s_i=(s_1,\ldots, s_n)$ for points $s_i\in \R$.) The dynamics
(\ref{mg}) can be viewed as an uncontrolled differential inclusion
$\dot x\in F(x)$ by taking $F(x)=\cup\{G(x,a): a\in A^n\}$. By
definition, the trajectories of (\ref{mg}) are those absolutely
continuous functions $y:[0,T]\to\R^n$ defined for some $T>0$ and
inputs $\alpha_i\in {\cal A}$ that satisfy $\dot y(t)\in
G(y(t),\alpha(t))$ for (Lebesgue) almost all (a.a.) $t\in [0,T]$.
The $D_i$'s can be interpreted as set valued state dependent
disturbance perturbations  acting on the individual components of
the locally Lipschitz dynamics $g=(g_1,g_2,\ldots, g_n)$. With
this interpretation, the values $t\mapsto \beta_i(t)\in D_i(y(t))$
assumed by the disturbances are unknown to the controller; one
only knows that each $\beta_i(t)$ takes {\em some value} in
$D_i(y(t))$ for
 each $t$. However, the mappings $g_i$ and $D_i$, the inputs
$t\mapsto \alpha_i(t)\in A$, and the current state $t\mapsto y(t)$
can be measured.  The dynamics (\ref{mg}) include the example from
the previous paragraph by taking $n=1$ and $g_1\equiv 1$. The
objective is to find sufficient conditions in terms of the $g_i$'s
and $D_i$'s under which all  trajectories of (\ref{mg}) starting
in a given closed set $S\subseteq\R^n$ remain in $S$, i.e., such
that $(G,S)$ is strongly invariant. Since the $D_i$'s are not
necessarily Lipschitz (or even continuous), the dynamics $G$ may
be discontinuous in the state, so the usual strong invariance
criteria
 for locally Lipschitz systems (cf.
\cite{AS03, CLSW98, V00}) do not apply.   Moreover, the dynamics
(\ref{mg}) are not in general tractable by the strong invariance
results from \cite{KMW04}, even if the $D_i$'s are singleton
valued.  For example, take $n=2$, $g_1\equiv g_2\equiv 1$,
$D_1(x)\equiv \{1\}$, $D_2(x)=\{1\}$ if $x_2\ge 0$, and
$D_2(x)=\{10\}$ if $x_2<0$.  In this case, if $f$ were  a feedback
realization for the $G$-trajectory $y(t)=(t,t)$ satisfying the
requirements of \cite{KMW04} (cf.  Section \ref{compare} for the
relevant definitions), then $(1,1)=f(t,(t,t))$ for all $t\ge 0$.
Since $f$ is continuous in the state variable, we can find
$\gamma>0$ such that $f(t,x)\in [1/2,3/2]^2$ for all $t\ge 0$ and
$||x||\le \gamma$.  In particular if $t\ge 0$, $||x||\le \gamma$,
and $x_2<0$, then $f(t,x)\not\in {\rm cone}\{D_1(x)\times
D_2(x)\}=\{\lambda(1,10): \lambda\ge 0\}$, so the cone requirement
on $f$  from \cite{KMW04} cannot be satisfied.
  We
prove strong invariance results for (\ref{mg}) in Section
\ref{sec4}.

Strong invariance theorems are of great independent interest
because they have been applied in many areas of nonlinear analysis
and dynamical systems theory.  For some applications, it suffices
to have characterizations of {\em weak invariance}, which is the
less restrictive requirement that for each point $\bar x$
 in the constraint set $S$, there exists {\em at least
one} trajectory $y$ of the dynamics starting at $\bar x$ such that
$y(t)\in S$ for all $t\ge 0$. While weak invariance
characterizations have been shown under very general assumptions
on the dynamics (e.g., locally bounded convex values and closed
graph), the standard results on strong invariance generally
require locally Lipschitz multifunctions. Strong invariance
theorems have important applications in uniqueness and regularity
theory for solutions of Hamilton-Jacobi-Bellman equations,
stability theory, differential games, monotone systems in biology,
and elsewhere (cf. \cite{AS03, AS04a, AS04b, CLSW98, CS03, V00}).
On the other hand, it is well appreciated that many important
dynamics such as (\ref{mg})  may be non-Lipschitz or even
discontinuous  in the state and so are beyond the scope of the
usual strong invariance methods. Therefore, the development of
conditions guaranteeing strong invariance under less restrictive
assumptions on the dynamics is a problem that is of considerable
ongoing research interest.

This motivates the search for new infinitesimal characterizations
for strong invariance for non-Lipschitz differential inclusions
which is the focus of this note. (Here and in the sequel,
``non-Lipschitz'' means ``not necessarily locally Lipschitz in the
state variable''.)   Donchev, Rios, and Wolenski \cite{DRW04a,
RW03} recently proved strong invariance characterizations under
the somewhat less restrictive structural assumption of one sided
Lipschitzness. A completely different approach was pursued by
Krastanov, Malisoff, and Wolenski \cite{KMW04} who gave a new
Hamiltonian sufficient condition for strong invariance for a class
of feedback realizable differential inclusions (see Section
\ref{compare} for the relevant definitions). The results in
\cite{KMW04}  do not require any of the usual structural
assumptions on the dynamics that are generally needed in strong
invariant systems theory, and therefore can be applied to a more
general class of systems. However, \cite{KMW04}  requires a cone
condition on the feedback realizations that is not in general
satisfied for the dynamics we consider below (see Section
\ref{compare} for more discussions on \cite{KMW04}).

In this note, we provide a nontrivial extension of  \cite{KMW04}
by proving  strong invariance under an alternative feedback
realizability condition that can in particular be satisfied by the
general non-Lipschitz system (\ref{mg}) and other examples that
are not tractable by  known strong invariance results. Our
condition is based  on Malisoff's coercive upper envelope approach
from \cite{M01} and Sussmann's unique limiting condition from
\cite{S03}.  We express our invariance results in terms of tangent
cones and Hamiltonian inequalities. Our approach has the
additional advantage over \cite{KMW04} that it is preserved under
``stacking'' in the following sense: If two dynamics $\dot y_1\in
G_1(y_1,\alpha_1)$ and  $\dot  y_2\in G_2(y_2,\alpha_2)$ satisfy
our feedback condition, then so does the ``stacked'' dynamic
$(\dot y_1,\dot y_2)\in G_1(y_1,\alpha_1)\times
G_2(y_2,\alpha_2)$.  The realizability condition in \cite{KMW04}
is not preserved under ``stacking''. Starting from  this
``stacking'' property, our results lead to new infinitesimal
monotonicity characterizations (see Section \ref{sec5} below).
Moreover, our results can still be applied even when the
constraint set $S$ is not necessarily a {\em closed} subset of
$\R^n$ (cf. Section \ref{relclosed} below).  While our main
theorem can be shown
 using  Zorn's Lemma, the proof we give below is constructive, and in
 particular leads to a
 new approach to building viable
trajectories for Carath\'eodory dynamics; see Remark \ref{franko}.

This note is organized as follows.  In Section \ref{sec2}, we
state our feedback realization property precisely and illustrate
its applicability  to  dynamics such as (\ref{mg}) that are not
tractable by the known strong invariance theory. In Section
\ref{sec3}, we state our new necessary and sufficient conditions
for strong invariance and we explain in detail how our conditions
improve on the known results. We prove our main results in Section
\ref{sec4}.  Section \ref{sec5} shows how our invariance
characterizations lead to
 new nonsmooth monotonicity characterizations for
systems with set valued state dependent non-Lipschitz
disturbances.
  We review the relevant arguments from \cite{KMW04} and
  background from nonsmooth
analysis in the appendices.

\section{Feedback Realizability Hypothesis and Examples}
\label{sec2}

 Our main object
 of study  is an autonomous
differential inclusion $\dot x\in F(x)$.  This includes measurable
controlled differential inclusions $\dot x\in G(x,a)$ with control
constraints $a\in A$
 by
taking $F(x):=G(x,A)$.  By a {\em trajectory} of $\dot x\in F(x)$
on an interval $[0,T]$ starting at a point $\bar x\in \R^n$, we
mean an absolutely continuous function $y :[0,T]\to \R^n$ for
which $y (0)=\bar x$ and $\dot y (t)\in F(y (t))$ for (Lebesgue)
almost all (a.a.) $t\in [0,T]$. We let ${\rm Traj}_T(F,\bar x)$
denote the set of all trajectories $y :[0,T]\to \R^n$ for $F$
starting at $\bar x$ for each $T>0$, and ${\rm Traj}(F,\bar
x):=\cup_{T \ge 0}{\rm Traj}_T(F,\bar x)$ and ${\rm
Traj}(F):=\cup_{\bar x\in \R^n} {\rm Traj}(F,\bar x)$. For a
control system  $\dot x\in G(x,A)$, we let ${\rm Traj}_T(G,\alpha,
\bar x)$ denote the set of all trajectories $y :[0,T]\to \R^n$ for
$G(\cdot, \alpha(t))$ starting at $\bar x$ for each $T>0$ and
$\alpha\in {\cal A}:=\{\, {\rm measurable}\;  \alpha:[0,T]\to
A\}$. In that case, we also set ${\rm Traj}(G,\alpha, \bar
x)=\cup_{T\ge 0}{\rm Traj}_T(G,\alpha, \bar x)$ and ${\rm
Traj}(G,\bar x )=\cup_{\alpha \in {\cal A}}{\rm Traj}(G,\alpha,
\bar x)$. For $n>1$ (respectively, $n=1$), a mapping $G$ defined
on a Borel subset of $\R^n$ is said to be {\em measurable}
provided $G$  is Borel (respectively, Lebesgue) measurable.

A multifunction $F:\R^n\rightrightarrows\R^n$ is said to have {\em
linear growth} provided there exist positive constants $c_1$ and
$c_2$ such that $||v||\le c_1+c_2||x||$ for all $v\in F(x)$ and
$x\in \R^n$, where $||\cdot||$ denotes the Euclidean supremum
norm. For any interval $I$, a function $f:I\times \R^n\to\R^n$ is
said to have {\em linear growth (on $I$)} provided $x\mapsto
F(x):=\{f(t,x): t\in I\}$ has linear growth. For any subsets
$C,P\subseteq \R^n$ and any constant $\eta\in \R$, we set $C+\eta
P:=\{c+\eta p: c\in C, p\in P\}$.  Also, ${\cal B}_n:=\{x\in \R^n:
\Vert x\Vert \le 1\}$.   A mapping $F:\R^n\rightrightarrows \R^n$
is said to be {\em lower semicontinuous} provided  for each $x\in
\R^n$ and $\eps>0$, there exists $\delta>0$ such that
$F(x')+\eps{\cal B}_n \supseteq F(x)$ for all $x'\in x+\delta
{\cal B}_n$; it is said to be {\em closed} (respectively,  {\em
compact}, {\em convex}, {\em nonempty}) {\em valued} provided
$F(x)$ is closed (respectively, compact, convex, nonempty) for
each $x\in\R^n$. We say that $F$ is {\em locally bounded} provided
$F(\eta{\cal B}_n)$ is bounded for each $\eta>0$. Throughout this
paper, we assume that all our mappings from $\R^n$ are nonempty
valued. Also, ${\rm int}(C)$ (respectively, ${\rm bd}(C)$) denotes
the interior (respectively, boundary) of any subset $C$ of a
Euclidean space.  The $i$th component of a mapping $F$ into $\R^n$
is denoted by $F_i$ for $i=1,2,\ldots, n$.

 A continuous
function $\omega:[0,\infty)\to[0,\infty)$ is called a {\em
modulus} provided it is nondecreasing  with $\omega(0)=0$. For
each $T\ge 0$,  we let $\C[0,T]$ denote the set of all functions
$f:[0,T]\times \R^{n}\to\R^{n}$ that satisfy

\begin{enumerate}\item[]\begin{itemize}
\item[$(C_1)$\ \ \ ]  For each $x\in\R^{n}$, the map $t\mapsto
f(t,x)$ is measurable;  \item[$(C_2)$\ \ \ ]  For each compact set
$K\subseteq\R^{n}$, there exists a modulus $\omega_{f,K}$ such
that for all $t\in[0,T]$ and $x_1,\,x_2\in K$, $\Vert
f(t,x_1)-f(t,x_2)\Vert \le \omega_{f,K}(\Vert x_1-x_2\Vert )$; and
 \item[$(C_3)$\ \ \ ] $f$ has linear growth on [0,T].
\end{itemize}\end{enumerate} This agrees with the definition of  $\C[0,T]$
in \cite{KMW04}. Our main hypothesis is that each  $y \in {\rm
Traj}(F)$ is also the unique (generalized) solution of an
appropriate initial value problem  $\dot x=f(t,x)$, $x(0)=y(0)$
for a feedback realization $f\in \C[0,T]$.  However, we allow a
more general class of feedback realizations than is allowed in
\cite{KMW04}.  We present our feedback realization hypothesis
next.

We need the following
additional definitions.  We let ${\rm co}$ (resp., $\overline {\rm
co}$) denote the convex hull (resp., closed convex hull).
 For each subset $P$ in Euclidean space,
we set ${\rm cone}\, \{P\}=\cup\{\eta P: \eta\ge 0\}$ (written
${\rm cone}\{p\}$ when $P=\{p\}$ is singleton)
 and we define   the  projections ${\rm pr}_i\{P\}=\{p_i:
 \exists p=(p_1, p_2, \ldots,
p_n)\in P\}$, so $F_i(x)={\rm pr}_i\{F(x)\}$ for each $i$ and
mapping $F$ into $\R^n$.
  We  define the {\em component
cone (ccone)} by ${\rm ccone}\{P\}=\{q\in\R^n: \exists p\in P \;
{\rm s.t.}\;  q_i\in {\rm cone}\{p_i\} \; \forall i\}$. When
$P=\{p\}$ is singleton, we write this as ${\rm ccone}\{p\}$.
Notice that $v\in {\rm ccone}\{P\}$ is a less restrictive
condition than $v\in {\rm cone}\{P\}$. For example, $(1,-1)\in
{\rm ccone}\{(2,-1)\}\setminus {\rm cone}\{(2,-1)\}$. Note that if
$P_1,P_2\subseteq \R^n$ and $v_i\in {\rm ccone}\{P_i\}$ for $i=1$
and $2$, then $(v_1,v_2)\in {\rm ccone}\{P_1\times P_2\}\subseteq
\R^{2n}$. Given $F:\R^n\rightrightarrows \R^n$, $\bar x\in \R^n$,
and $T>0$, we set
\[\C'_F([0,T], \bar x):=\{f\in \C[0,T]: \exists \gamma>0 {\rm\
s.t.\ } f(t,x)\in {\rm ccone} \{F(x)\} \text{ for a.a.
}t\in[0,\gamma) \; \text{and all} \; x\in \bar x+\gamma {\cal
B}_n\}.\] Given $W\subseteq\R^n$, we call a mapping
$F:\R^n\rightrightarrows\R^n$  {\em $W$-weakly zeroing} provided:
$[x\in W\, \;  {\rm and}\, \;  p\in F(x)]$ $\Rightarrow$
$\{(p_1,0,0,\ldots,0), (0,p_2,0,\ldots, 0),\ldots, (0,0,\ldots, 0,
p_n)\}\subseteq F(x)$. We simply say that $F$ is {\em weakly
zeroing} provided it is $\R^n$-weakly zeroing.  Any $F$ that is
the product of one dimensional multi-functions $F_i$ satisfying
$0\in F_i(x)$ for all $i$ and $x$ is weakly zeroing, but weak
zeroing  is more general since it allows images $F(x)$ such as
$\{p\in \R^2: p_1\ge 0, p_2\ge 0, p_2\le 1-p_1\}$. We assume the
following condition:

\begin{itemize}\item[]\begin{itemize}
\item[$(U')$\ \ \ ] For each $\bar x$, $T\ge 0$, and $y\in {\rm
Traj}_T(F,\bar x)$, there exists $f\in \C'_F([0,T], \bar x)$ for
which $y$ is the unique solution   of  $\dot y(t)=f(t,y(t))$ on
$[0,T]$  starting at $\bar x$.  Also, $F$ is weakly zeroing.
\end{itemize}\end{itemize}

When the conditions in $(U')$ hold, we call $f$ a {\em feedback
realization} for the trajectory $y$. The prime notation signifies
that our condition is a variant of the realizability condition
$(U)$ from \cite{KMW04}; see  Section \ref{compare} for a
comparison of conditions $(U)$ and $(U')$.  Since  $0\in F(x)$ for
all $x\in \R^n$, our dynamics are {\em weakly} (but not
necessarily strongly) invariant for any constraint set $S$. In our
tangential characterizations of strong invariance, our weak
zeroing requirement can be relaxed to the requirement that $F$ be
${\rm bd}(S)$-weakly zeroing for the closed constraint set $S$
(see Section \ref{suff}).  Notice that feedback realizations $f$
can depend on the initial value $\bar x$ and the trajectory $y$.
For example, if $F$ is compact-convex valued and Lipschitz and
weakly zeroing, then (see \cite{KMW04}) we can satisfy $(U')$
using $f(t,x):= {\rm proj}_{F(x)}\{\dot y(t)\}$ for any $y\in {\rm
Traj}(F)$, where ${\rm proj}_{F(x)}\{q\}$ denotes the (unique)
closest point in $F(x)$ to $q$.
 The following {\em non-Lipschitz} examples  show how Condition $(U')$
 can be satisfied without
 necessarily having $f(t,x)\in F(x)$ for all $x$ and a.a. $t$.
\begin{example}
\label{ex0.5} Choose $n=1$, $F(0)=[-1,+1]$, and $F(x)=\{-{\rm
sign}(x),0\}$ for $x\ne 0$. To verify Condition $(U')$, let $T>0$,
$\bar x\in \R$, and $y\in {\rm Traj}_T(F,\bar x)$ be given. Note
that $(F,\{0\})$ is strongly invariant. Therefore, either (i)
$y(t)$ starts at some $\bar x\ne 0$ and then moves toward $0$ at
unit or zero speed or (ii) $y(t)\equiv 0$. If $\bar x\ne 0$, then
$(U')$ is satisfied using  $f(t,x)\equiv -{\rm sign}(\bar x)
1_D(t)$, where $1_D(t)$ is defined to be   $0$ when $\dot y(t)=0$
or $y(t)=0$, and $1$ otherwise. Then $f(t,x)\in
 {\rm cone}\{F(x)\}$ for all $t\in [0,T]$ and
$x\in \bar x+(|\bar x|/2){\cal B}_1$.   If instead $\bar x =0$,
then we choose  $f(t,x)\equiv 0\in {\rm cone}\{F(x)\}$ for all
$t\in [0,T]$ and $x\in \R$.
\end{example}

\begin{example}
\label{ExhibitA} Consider the following dynamics evolving on
$\R^n$:
\begin{equation}
\label{exdyn} \dot x\in F(x)=\prod_{i=1}^n L_i(x)D_i(x),
\end{equation}
where we assume the following for each $i$:
\begin{itemize}\item[]\begin{itemize}
\item[$(H_1)$\ \ \ ] $x\mapsto L_i(x)\subseteq\R$ is locally
Lipschitz and compact valued. \item[$(H_2)$\ \ \ ]  $x\mapsto
D_i(x)\subseteq\R$ is measurable, locally bounded, and closed
valued. \item[$(H_3)$\ \ \ ] $0\in D_i(x)$ for all $x\in \R^n$,
and $x\mapsto {\rm cone}\{D_i(x)\}$ is constant.
\end{itemize}\end{itemize}
where $L_i(x)D_i(x):=\{\lambda \delta: \lambda\in L_i(x),
\delta\in D_i(x)\}$.
 Note that under $(H_3)$, the only possible constant
values for ${\rm cone}\{D_i(x)\}$ are $\{0\}$, $\R$, $[0,\infty)$,
or $(-\infty,0]$.  In particular, for each $x,y\in \R^n$ and
$\beta\in D_i(x)$, we have $\beta\in {\rm cone}\{D_i(y)\}$.
 These hypotheses allow the systems
(\ref{mg}) from the introduction by taking $L_i(x)\equiv g_i(x,A)$
for all $i$.   To check that (\ref{exdyn}) satisfies   $(U')$, let
$\bar x\in \R^n$, $T>0$, $y\in {\rm Traj}_T(F, \bar x)$, and
$\phi:\R^n\to [0,1]$ be any smooth (i.e., $C^\infty$)  function
that is identically one on  $y$ and compactly supported. Then $y$
is also a trajectory of
\begin{equation}
\label{s1} G(x):=\prod_{i=1}^n \phi(x)  \, {\rm co} \{L_i(x)\}\,
D_i(x).
\end{equation}
Since each component $x\mapsto \phi(x) {\rm co}\{L_i(x)\}$ is
Lipschitz and compact  and convex valued, standard parametrization
results (cf. \cite[Section 9.7]{AF90}) allow us to rewrite
(\ref{s1}) as $G(x)=\prod_{i=1}^n g_i(x, {\cal B}_1) D_i(x)$ for
some Lipschitz function $g=(g_1,g_2,\ldots, g_n)$ where
$g_i:\R^n\times {\cal B}_1\to \R$ for all $i$.  Since $y\in {\rm
Traj}_T(G, \bar x)$, there are measurable functions
$\alpha_i:[0,T]\to {\cal B}_1$ and $\beta_i:[0,T]\to \R$, with
$\beta_i(t)\in D_i(y(t))$ for a.a. $t\in [0,T]$, such that $y$ is
the unique solution of
\begin{equation}\label{uni}f(t,x):=\prod_{i=1}^n g_i(x, \alpha_i(t))
\beta_i(t)\end{equation}  on $[0,T]$ starting at $\bar x$. This
follows from the (Generalized) Filippov Selection Theorem (see
\cite[p.72]{V00}).
  The local boundedness
requirement in $(H_3)$ guarantees that the $\beta_i$'s are
(essentially) bounded. Using the fact that ${\rm cone}\{{\rm
co}D\}={\rm cone}\{D\}$ for all subsets $D\subseteq\R$, one can
verify  $(U')$ using the choice (\ref{uni}) of $f$.
\end{example}
\begin{example}
\label{stacked} We next illustrate how dynamics satisfying $(U')$
can be ``stacked'' to build new dynamics that again satisfy
$(U')$. Assume two dynamics $F_1,F_2:\R^n\rightrightarrows\R^n$
satisfy condition  $(U')$. Consider the ``stacked'' dynamics
$F:\R^{2n}\rightrightarrows \R^{2n}$ defined by
$F(x)=F_1(x_1)\times F_2(x_2)$.  We claim that $F$ satisfies
$(U')$ in dimension $2n$. To see why, let $T>0$, $\bar x=(\bar
x_1, \bar x_2)\in \R^n\times \R^n$, and $y\in {\rm Traj}_T(F,\bar
x)$ be given. Write $y=(y_1,y_2):[0,T]\to \R^n\times \R^n$. Since
$y_i\in {\rm Traj}_T(F_i, \bar x_i)$ for $1\le i\le 2$, we can
find $f_i\in C'_{F_i}([0,T],\bar x_i)$ such that $y_i$ is the
unique solution of the initial value problem $\dot
x_i=f_i(t,x_i)$, $x_i(0)=\bar x_i$ on $[0,T]$ for $1\le i\le 2$.
The  requirement $(U')$ is then satisfied for $y$ using the
feedback realization $f(t,x):=(f_1(t,x_1), f_2(t,x_2))$.  On the
other hand, since we do not in general have $f(t,x)\in {\rm
cone}\{F(x)\}$, the stacked dynamics may not satisfy the more
restrictive requirement $(U)$ from \cite{KMW04} (see Section
\ref{compare}).

\end{example}

The preceding examples  play an important role when we apply our
results to monotone control systems in Section \ref{sec5}.  In the
next two sections, we state and prove our main strong invariance
results. For the relevant background from nonsmooth analysis, see
Appendix \ref{nonsmo}.

\section{Statement and Discussion of Main Results}
\label{sec3}
\subsection{Sufficient Conditions for Strong Invariance}
\label{suff}
  Let
$H_F:\R^n\times\R^n\to [-\infty, +\infty]$ denote the {\em (upper)
Hamiltonian} for our dynamics $F$ defined by
$H_F(x,d):=\sup\{\langle v,d \rangle: v\in F(x)\}$. For any subset
$D\subseteq\R^n$, we write $H_F(x,D)\le 0$ to mean that
$H_F(x,d)\le 0$ for all $d\in D$.  By definition, this inequality
holds vacuously if $D=\emptyset$.

\begin{theorem}
\label{main} Let $F:\R^n\rightrightarrows \R^n$  satisfy $(U')$
and  $\Psi:\R^n\to \R$ be lower semicontinuous. If there exists an
open set ${\cal U}\subseteq \R^n$ containing ${\cal S}:=\{x\in
\R^n: \Psi(x)\le 0\}$  for which $H_F(x,\partial_P\Psi(x))\le 0$
for all $x\in {\cal U}$, then $(F,{\cal S})$ is strongly
invariant.
\end{theorem}
The statement of Theorem \ref{main} is the same as the main
theorem  in \cite{KMW04} except Theorem \ref{main}  replaces the
realization assumption $(U)$ from \cite{KMW04} with  $(U')$ (see
Section \ref{compare} for a comparison of these two assumptions).
{}We prove Theorem \ref{main} below by constructing appropriate
Euler polygonal arcs;  see also Remark \ref{franko} for an
alternative nonconstructive  proof based on Zorn's Lemma.
Condition $(U')$ allows us to apply our results to a broad class
of examples
 that are not tractable
by \cite{KMW04} or other known strong invariance results (see
Section \ref{sec2}).
 Note that we require the Hamiltonian
inequality in a {\em neighborhood} ${\cal U}$ of ${\cal S}$.  The
result is not true in general if the Hamiltonian condition is
placed only on ${\cal S}$. For example, take $n=1$, $\Psi(x)=x^2$,
and $F(x)\equiv \{0,1\}$. Then ${\cal S}=\{0\}$ and
$H_F(0,\partial_P\Psi(0))=0$, but in this case $(F, {\cal S})$ is
not strongly invariant.  The non-Lipschitz system in Example
\ref{ex0.5} is covered by Theorem \ref{main} if we choose the
verification function $\Psi(x)=x^2$. In that case, the Hamiltonian
condition reads $H_F(x,\Psi'(x))\equiv  0$ for all $x\in \R$, so
our sufficient condition for strong invariance is satisfied.

Theorem \ref{main} contains the usual sufficient condition for
strong invariance for an arbitrary closed set $S\subseteq\R^n$ by
letting $\Psi$ be the characteristic function $\chi_S$ of $S$;
that is, $\chi_S(x)=0$ if $x\in S$ and is $1$ otherwise. Then
$\partial_P\Psi(x)=\{0\}$ for all $x\not\in {\rm bd}(S)$, and
$\partial_P\Psi(x)=N^P_S(x)$ for all $x\in {\rm bd}(S)$.
 This implies:
\begin{corollary}
\label{cor1} Let $F:\R^n\rightrightarrows \R^n$  satisfy $(U')$
and $S\subseteq \R^n$ be closed.  If $H_F(x,N^P_S(x))\le 0$ for
all $x\in {\rm bd}(S)$, then $(F,S)$ is strongly invariant.
\end{corollary}

In fact, Corollary \ref{cor1} is   the special case of Theorem
\ref{main} using ${\cal U}=\R^n$, $\Psi=\chi_S$, and ${\cal S}=S$.
 For the proof of Theorem \ref{main}, see Section \ref{sec4}.
The proof shows that Theorem \ref{main} remains true even if the
weakly zeroing requirement on $F$ is relaxed to:  There exists an
open set ${\cal G}$ containing ${\rm bd}(S)$ such that $F$ is
${\cal G}$-weakly zeroing. Moreover, Corollary \ref{cor1} remains
true even if this requirement is relaxed to requiring that $F$ be
${\rm bd}(S$)-weakly zeroing.

\subsection{Strong Invariance Characterizations}\label{relclosed}
 The converse of
Corollary \ref{cor1} does not hold, as illustrated by  Example
\ref{ex0.5}.  There $(F,\{0\})$ is strongly invariant but
$N^P_{\{0\}}(0)=\R$, so the Hamiltonian condition in the corollary
is not satisfied. This means that the converse of Theorem
\ref{main} does not hold. On the other hand, it is desirable to
have  a sufficient {\em and necessary} condition for strong
invariance, that is, a strong invariance characterization. One
would also hope to have such a characterization  in the more
general situation where a non-Lipschitz dynamic $F$ evolves on an
arbitrary open subset ${\cal O}\subseteq \R^n$ and the constraint
set $S\subseteq\R^n$ is a (relatively) closed subset of ${\cal O}$
but is not necessarily a closed subset of $\R^n$; i.e.,  $S={\cal
F}\cap {\cal O}$ for some closed subset ${\cal F}\subseteq \R^n$.
Characterizations of this kind play an important role in the
analysis of monotone control systems (see Section \ref{sec5} below
or \cite{AS03}).

 In this subsection, we  provide such a
characterization in terms of the Bouligand tangent cone $T^B_S$
(see Appendix \ref{nonsmo} for the relevant definitions). Strong
invariance of $(F,S)$ under relatively closed state constraints
was treated in \cite{AS03}. There it was assumed that $F$ is
locally Lipschitz.
  Here we consider  a more general dynamic
$F:{\cal O}\times A\rightrightarrows \R^n$  of the form
\begin{equation}
\label{relaform}
 \dot x\in F(x,a)=\prod_{i=1}^n g_i(x,a)D_i(x), \; \; x\in {\cal
 O},\; \; a\in A
\end{equation}
where $A$ is a compact subset of a Euclidean space, each
$g_i:{\cal O}\times A\to \R$ is locally Lipschitz, and
\begin{itemize}\item[]\begin{itemize}
\item[$(H_4)$\ \ \ ]  Each $D_i$ is convex valued, lower
semicontinuous, and satisfies $(H_2)$-$(H_3)$.
\end{itemize}\end{itemize}

 \noindent Our assumptions imply that ${\cal O}\ni
x\mapsto F(x,A)$
 satisfies $(U')$ (see
  Example \ref{ExhibitA}).
 We prove:

\bt{relathm} Let ${\cal O}\subseteq \R^n$ be open and  $S$ be a
(relatively) closed subset of ${\cal O}$.  Let $F:{\cal O}\times
A\rightrightarrows \R^n$ be as above.
 Then $(F,S)$ is
strongly invariant  if and only if $F(x,A)\subseteq T^B_S(x)$ for
all $x\in S$.\et

We defer the proof to Section \ref{sec4} which also shows that the
requirement that $0\in D_i(x)$ for all $i$ and $x$ can be relaxed
to: There exist compact sets  $N_1\subseteq N_2\subseteq
\ldots\subseteq N_k\subseteq\ldots$  such that
\begin{itemize}
\item[(i)]  ${\cal O}=\cup_k N_k$ and \item[(ii)] $0\in D_i(x)$
for all $i$, all $x\in {\rm bd}(S\cap N_k)$, and all $k\in
\N$,\end{itemize} which imply that the pairs $(F, S\cap N_k)$ are
weakly invariant.

\begin{remark}
\label{convex} The tangential condition from Theorem \ref{relathm}
remains {\em sufficient} for strong invariance of $(F,S)$ if
$(H_4)$ is relaxed to requiring $x\mapsto F(x,A)$ to satisfy
$(U')$, by the same proof.
 This gives a sufficient condition for strong invariance
for the dynamics (\ref{relaform}) with relatively closed state
constraints if we merely assume $(H_2)$-$(H_3)$ without
necessarily having convex valued or lower semicontinuous dynamics.
 On the other hand, Example \ref{ex0.5} illustrates
that the tangential condition is not {\em necessary} for strong
invariance without assuming $(H_4)$ even if ${\cal O}=\R^n$. In
that example, there are no controls $a$, and $(F, \{0\})$ is
strongly invariant but $[-1,+1]=F(0)\not\subseteq
T^B_{\scriptscriptstyle \{0\}}=\{0\}$. This does not contradict
Theorem \ref{relathm} since  $F$ is not lower semicontinuous, so
$(H_4)$ is not satisfied.\end{remark}

\subsection{Relationship to Known Strong Invariance Results}
\label{compare} Theorem \ref{main}  applies to a more general
class of multifunctions than
 the known strong invariance results because it does
not require the usual Lipschitz, one-sided Lipschitz,  or other
structural assumptions on $F$, nor does it require the more
restrictive feedback realizability condition from \cite{KMW04}.
The papers \cite{C75, CLR97} provide strong invariance results for
locally Lipschitz dynamics (see also \cite[Chapter 4]{CLSW98}).
For locally Lipschitz $F$, Clarke \cite{C75} showed that strong
invariance  of $(F,S)$ is equivalent to $F(x)\subseteq T^C_S(x)$
for all $x\in S$, where $T^C_S$ denotes the Clarke tangent cone
(cf. \cite{CLSW98} or Appendix \ref{nonsmo} for the relevant
definitions).
 See \cite{CLR97} for Hilbert space versions and
\cite{V00} for other strong invariance results for Lipschitz
dynamics and nonautonomous versions.  For strong invariance
characterizations under somewhat more general structural
conditions on $F$ (e.g., dissipativity and one-sided
Lipschitzness), see \cite{D02, DRW04a,  RW03}.

On the other hand,  Theorem \ref{main} does not make any such
structural assumptions on $F$ and allows non-Lipschitz dynamics
such as (\ref{exdyn}) that are intractable by the standard strong
invariance results. In \cite{KMW04}, strong invariance was shown
assuming the following on the multifunction $F$:

\begin{itemize}\item[]\begin{itemize}
\item[$(U)$\ \ \ ] For each $\bar x\in \R^n$, $T\ge 0$, and $y\in
{\rm Traj}_T(F,\bar x)$, there exists $f\in \C_F([0,T], \bar x)$
for which $y$ is the unique solution   of the initial value
problem $\dot y(t)=f(t,y(t))$, $y(0)=\bar x$ on $[0,T]$.
\end{itemize}\end{itemize}
where $\C_F([0,T], \bar x)= \left\{f\in \C[0,T]: \exists \gamma>0
{\rm\ s.t.\ } f(t,x)\in {\rm cone}\, \{F(x)\}\, \forall x\in \bar
x+ \gamma {\cal B}_n \; \& {\rm \ a.a.\  }\, t\in[0,T]\right\}$
and $\C[0,T]$ is as defined in Section \ref{sec2} above.  {\em An
important difference between $(U)$ and $(U')$ is that $(U)$
requires $f(t,x)$  to be locally in the cone of $F(x)$ rather than
the larger ccone so Condition $(U)$ from \cite{KMW04} is a more
restrictive requirement than our Condition $(U')$.} In other
words, our condition $f(t,x)\in {\rm ccone}\{F(x)\}$ from $(U')$
means there are weights $\omega_1,\omega_2,\ldots, \omega_n\in
[0,\infty)$ (possibly depending on $t$ and $x$) such that
$f_i(t,x)\in \omega_i F_i(x)$ for all $i$ while $(U)$ makes the
further restriction that $\omega_1=\omega_2=\ldots=\omega_n$ (see
Section \ref{sec1} for an example where $(U')$ holds but these
weights cannot be chosen to be equal).
 The main results of \cite{KMW04} are the same as our Theorem
\ref{main} and Corollary \ref{cor1} except with $(U')$ replaced by
$(U)$.

Condition $(U)$ has an advantage because it does not require
zeroing, but the argument in the next section shows that Corollary
\ref{cor1} remains true even if this requirement is relaxed to
requiring that $F$ be ${\rm bd}(S)$-weakly zeroing.
 Notice
however that if $F_1,F_2:\R^n\rightrightarrows \R^n$ are two
dynamics that satisfy $(U)$, then it is not in general the case
that the ``stacked'' dynamic $F(x):=F_1(x_1)\times F_2(x_2)$
satisfies $(U)$. In particular, if $\bar x=(\bar x_1, \bar x_2)\in
\R^n\times \R^n$, and if also $f_1,f_2\in \C[0,T]$ and $\gamma>0$
are such that $f_j(t,x_j)\in {\rm cone}\{F_j(x_j)\}$ for all
$x_j\in \bar x_j+\gamma {\cal B}_n$, a.a. $t\in [0,T]$, and
$j=1,2$, then it is {\em not} generally the case that
$(f_1(t,x_1), f_2(t,x_2))\in {\rm cone}\{F_1(x_1)\times
F_2(x_2)\}$, although we do have
\[(f_1(t,x_1), f_2(t,x_2))\in {\rm ccone}\{F_1(x_1)\times
F_2(x_2)\}\; \;  {\rm for\ \ all }\; \;  x=(x_1,x_2)\in \bar
x+\gamma{\cal B}_{2n}\; \;  {\rm and\ \   a.a.}\;  \; t\in
[0,T],\] so $F$ satisfies $(U')$ in dimension $2n$. Therefore,
condition $(U')$ has the advantage that it is preserved under
``stacking''. Using this ``stacking'' property, we  prove
nonsmooth monotonicity characterizations that extend the
corresponding results from \cite{AS03} (see Section \ref{sec5}).
Theorem \ref{relathm} extends the strong invariance
characterization for locally Lipschitz controlled dynamics and
relatively closed state constraints from \cite[Theorem 4]{AS03} by
allowing non-Lipschitz disturbances $D_i$. In particular, Theorem
\ref{relathm} gives a tangential characterization for strong
invariance for non-Lipschitz dynamics defined on all of ${\cal
O}=\R^n$; our results are new even for this particular case.

\section{Proof of Main Results}
\label{sec4}

\subsection{Proof of Theorem \ref{main}}
\label{4p1}
 Fix $T>0$, $\bar x\in {\rm bd}({\cal S})$, $\eps>0$,
and $f\in C'_F([0,T],\bar x)$.   By definition, we can then find
$\gamma\in (0,1)$ such that $f_i(t,x)\in {\rm cone}\{F_i(x)\}$ for
all $i$, a.a. $t\in [0,\gamma)$, and all $x\in \bar x+\gamma{\cal
B}_n$. We define the mollification $f_\eps$ by (\ref{ftild}) in
Appendix \ref{nonsmo} below and  let $f_{\eps,i}$ denote its $i$th
component for $i=1,2,\ldots, n$. By reducing $\gamma$ as necessary
without relabelling, we can assume that $\bar x+\gamma{\cal
B}_n\subseteq {\cal U}$.  By also reducing $T>0$, we also assume
$T\in (0,\gamma)$. For each $i\in \{1,2,\ldots, n\}$, $t\in
[0,T]$, $x\in \R^n$, and $r>0$, we then set
\[
G^\eps_{f,i}[t,x,r]=\prod_{j=1}^nR^i_j, \; \; \;
G^\eps_f[t,x,r]=\prod_{j=1}^n G^\eps_{f_j}[t,x,r],\; \; \;
g_f[t,x,r]:=1+\sup\left\{||p||: p\in G^\eps_{f}[t,x,r]\right\},
\]
where $R^i_j\equiv \{0\}$ for all $j\ne i$ and
\[R^i_i=G^\eps_{f_i}[t,x,r]:=\overline{\rm co}\{f_{\eps,i}(t,y):
||y-x||\le 1/r\}.\] It follows that  $G^\eps_{f,i}[t,x,r]\subseteq
g_f[t,x,r]{\cal B}_n$ everywhere for each $i$.   By our linear
growth assumption $(C_3)$ from Section \ref{sec2}, the sets
$G^\eps_{f,i}$ are compact. The following consequence of the
Clarke-Ledyaev Mean Value Inequality (Theorem \ref{MVI} from
Appendix \ref{nonsmo}) extends Claim 4.1 from \cite{KMW04} to our
more general feedback realizations $f$:

\bcl{claim1} If $x\in \R^n$, $t\ge 0$, $k\in \N$, and $h>0$ are
such that
\begin{equation}\label{crucial}
0 < h \le \frac{1}{32k\, g_f[t,x,k/n]}\quad\text{and}\quad x+h
g_f[t,x,k/n]{\cal B}_n \subseteq \bar x+\frac{2\gamma}{3}{\cal
B}_n,
\end{equation}
then
\begin{equation}\label{crucialConclusion}
\Psi(x+hv)\le \Psi(x)+\frac{h}{k}
\end{equation}
holds for some $v\in G_{f}^\eps[t,x,k/n]$. \ecl

\begin{proof}
Set $Y_1= x+hG^\eps_{f,1}[t,x,16k]$ which is compact and convex.
We first find $v_1\in G^\eps_{f_1}[t,x,16k]$ such that
\begin{equation}\label{numb1}
\Psi(x+h(v_1,0,\ldots,0))\le \Psi(x)+\frac{h}{nk}.
\end{equation}
Suppose no such $v_1$ exists.  Since $\Psi$ is lower
semicontinuous and $Y_1$ is compact, we would then have
\begin{equation}
\label{ddef} \delta :=\frac{h}{nk}< \min_{y\in
Y_1}\Psi(y)-\Psi(x).
\end{equation}
 Notice that
$g_f[t,x,16k]\le g_f[t,x,k/n]$.  By (\ref{crucial}), we
 can therefore find $\lambda\in (0,\frac{1}{32k})$ satisfying
\begin{equation}
 \label{keyy}
x+hg_f[t,x,16k]{\cal B}_n+\lambda{\cal B}_n \; \; \subseteq \; \;
x+hg_f[t,x,k/n]{\cal B}_n+\lambda{\cal B}_n \; \; \subseteq \; \;
\bar x+\gamma{\cal B}_n\; \; \subseteq\; \;  {\cal
U}.\end{equation}
 Next we apply Theorem
\ref{MVI} from Appendix \ref{nonsmo}  with the choices $Y=Y_1$ and
$\delta$ defined by (\ref{ddef}). It follows that there exist $z
\in [x,Y_1] + \lambda{\cal B}_{n}$ and $\zeta\in
\partial_P\Psi(z)$ for which
\begin{equation}\label{this}
\delta<\min_{y\in Y_1}\langle \zeta, y-x\rangle= \min_{v\in
G^\eps_{f,1}[t,x,16k]}\langle \zeta,hv\rangle,
\end{equation}
where $[x,Y_1]$ denotes the closed convex hull of $x$ and $Y_1$.
By (\ref{keyy}),  $z\in \bar x+\gamma{\cal B}_n\subseteq {\cal
U}$.  Also,   (\ref{crucial}) combined with the choice of
$\lambda$ gives $ \Vert z-x\Vert \; \le \; hg_f[t,x,16k]+\lambda\;
\le\; hg_f[t,x,k/n]+\lambda\; \le \; \frac{1}{16k}$. Therefore
\begin{equation}\label{np}f^o_{\eps,1}(t,z):=(f_{\eps,1}(t,z),0,0,\ldots, 0) \in
G^\eps_{f,1}[t,x,16k].\end{equation} Since $z\in \bar
x+\gamma{\cal B}_n$, we know that $f_1(s,z)\in {\rm
cone}\{F_1(z)\}$ for a.a. $s\in [0,T]$. Since $F$ is weakly
zeroing, this gives $(f_1(s,z),0,0,\ldots, 0)\in {\rm
cone}\{F(z)\}$ for a.a. $s\in [0,T]$.  Since $z\in {\cal U}$,  our
Hamiltonian hypothesis then gives
\begin{equation}\label{kee2}\langle \zeta, (f_1(s,z),0,0,\ldots, 0)\rangle \le 0\; \; \;
{\rm for\ \ a.a.}\; \;  s\in [0,T].\end{equation} Therefore,
 (\ref{this})-(\ref{kee2}) give the contradiction
\begin{equation}
\label{contra} 0\; \; <\; \; \delta\; \; \leq\; \; h\langle \zeta,
f^o_{\eps,1}(t,z)\rangle\; \; =\; \; h
\int_{\R}\eta_\eps(t-s)\langle \zeta, (f_1(s,z),0,0,\ldots,
0)\rangle ds\; \; \le\; \; 0,
\end{equation}
so there must exist $v_1\in G^\eps_{f_1}[t,x,16k]\subseteq
G^\eps_{f_1}[t,x,k/n]\subseteq g_f[t,x,k/n]{\cal B}_1$ satisfying
(\ref{numb1}). Assume $n>1$.

Set $v^o_1=(v_1,0,0,\ldots, 0)\in \R^n$.  Then $(x+hv^o_1)  +
hg_f[t,x,k/n](\{0\}\times {\cal B}_{n-1}) \subseteq x +
hg_f[t,x,k/n]{\cal B}_{n}$,  so (\ref{keyy}) gives $ (x+hv^o_1) +
hg_f[t,x,k/n](\{0\}\times {\cal B}_{n-1}) +  \lambda {\cal B}_n
\subseteq  \bar x+\gamma{\cal B}_n$.
 Moreover, \begin{equation}\label{tr}G^\eps_{f_2}[t,x+hv^o_1,16k]\subseteq
 G^\eps_{f_2}[t,x,k/n],\end{equation}  because
 if  $||y-(x+hv^o_1)||\le
\frac{1}{16k}$, then (\ref{crucial})  gives $||y-x||\; \le\;
\frac{1}{16k}+h|v_1|\; \le \; \frac{1}{16k}+\frac{1}{32k}\; \le\;
\frac{2}{16k}\; \le\; \frac{n}{k}$, which gives (\ref{tr}).
 We can therefore apply the preceding argument with $x$ replaced
by $x+hv^o_1$, and with $Y_1$ replaced by the new compact convex
set
\[Y_2\; \; :=\; \; (x+hv^o_1) +hG^\eps_{f,2}[t,x+hv^o_1,16k]\; \;
\subseteq\; \; (x+hv^o_1)+hg_f[t,x,k/n](\{0\}\times {\cal
B}_{n-1})
\] to find $z\in [x+hv^o_1,Y_2]+\lambda {\cal B}_n$ such that
$f^o_{\eps,2}(t,z):=(0, f_{\eps,2}(t,z),0,0,\ldots, 0) \in
G^\eps_{f,2}[t,x+hv^o_1,16k]$ and  $v_2\in
G^\eps_{f_2}[t,x+hv^o_1,16k]\subseteq G^\eps_{f_2}[t,x,k/n]$ such
that
\begin{equation}\label{numb2}
\Psi(x+hv^o_1+h(0,v_2,0,0,\ldots, 0))\le
\Psi(x+hv^o_1)+\frac{h}{nk}.
\end{equation}
Next we argue by induction.  Proceeding inductively, we  apply the
same argument but with $x$ replaced by $\tilde
x_{i-1}:=x+h(v_1,v_2,\ldots, v_{i-1},0,0,\ldots, 0)$ and with the
set $Y_1$ replaced by the new compact convex set $Y_i:=\tilde
x_{i-1}+hG^\eps_{f,i}[t, \tilde x_{i-1}, 16k]$. Since $v_r\in
G^\eps_{\scriptscriptstyle f_r}[t,x,k/n]$ for all $r<i$,
(\ref{crucial}) implies that $h|v_r|\le \frac{1}{32k}$ for all
$r<i$, so the proof of (\ref{tr}) shows  $G^\eps_{f,i}[t,\tilde
x_{i-1}, 16k]
 \subseteq G^\eps_{f,i}[t,x,k/n]$ and
\[
\begin{array}{lll}
[\tilde x_{i-1}, Y_i]+\lambda {\cal B}_n &\subseteq& \tilde
x_{i-1} +hg_f[t,x,k/n](\{0\}\times {\cal B}_{n-i+1})+\lambda {\cal
B}_n\\
& \subseteq& x+(hg_f[t,x,k/n]+\lambda){\cal B}_n\; \; \subseteq\;
\;  \bar x+\gamma {\cal B}_n\; \; \subseteq\; \; {\cal U} \; \;
({\rm by\ } (\ref{keyy})).
\end{array}\]
This allows us to find
 $v_i \in G^\eps_{f_i}[t,\tilde x_{i-1}, 16k]
 \subseteq G^\eps_{f_i}[t,x,k/n]$ such
that
\begin{equation}\label{gd}
\Psi(x+h(v_1,v_2,\ldots, v_i,0,0,\ldots, 0))\le \Psi(\tilde
x_{i-1})+\frac{h}{nk}.
\end{equation}
for $i=2,3,\ldots, n$. Choosing $v=(v_1,v_2,\ldots, v_n)\in
G^\eps_{f}[t,x,k/n]$, we obtain (\ref{crucialConclusion}) by
summing the inequalities in (\ref{numb1})  and (\ref{gd}) over
$i=2,\ldots, n$.
\end{proof}

Now set $D:= \bar x+\frac{\gamma}{2}{\cal B}_n\subseteq {\cal U}$.
Let $\omega_{f,K}$ be a modulus of continuity for $x\mapsto
f(t,x)$ on $K:=D+n{\cal B}_n$ for all $t\in [0,T]$ (see Condition
$(C_2)$). Then $\omega_{f,K}$ is also a modulus of continuity of
$K\ni x\mapsto f_\eps(t,x)$
for all $t\in [0,T]$ and $\eps>0$.  The following claim parallels
Claim 4.2 from \cite{KMW04}:

\bcl{claim2}

Let $(t,x,k)\in [0,T]\times D\times \N$ and $v\in
G^\eps_{f}[t,x,k/n]$. Then $\Vert v-f_\eps(t,x)\Vert \le
\omega_{f,K}(n/k)+1/k$.\ecl

\bpr Using the Carath\'eodory Theorem  and the definition of
$G^\eps_f[t,x,k/n]$, we can write
\[
v=\Delta+\prod_{i=1}^n\left[ \sum_{j=1}^{2} \alpha_{i,j}
f_{\eps,i}(t,y_{i,j})\right],
\]
where \[\Vert y_{i,j}-x\Vert \le n/k \; \forall i,j;\; \; \;
\alpha_{i,j}\in [0,1]\; \; \forall i,j;\; \; \; \sum_{j=1}^{2}
\alpha_{i,j}=1\; \; \forall i; \; \; \;  {\rm and}\; \; \;
||\Delta||\le 1/k.\] In particular, $x\in K$ and $y_{i,j}\in K$
for all $i$ and $j$. This gives
\[
\displaystyle
\begin{array}{lll}
|v_i-f_{\eps,i}(t,x)| \;  \le \; \left\vert
v_i-\displaystyle\sum_{j=1}^{2} \alpha_{i,j}
f_{\eps,i}(t,y_{i,j})\right\vert\;  +\;
\left\vert\displaystyle\sum_{j=1}^{2} \alpha_{i,j}
\{f_{\eps,i}(t,y_{i,j})-f_{\eps,i}(t,x)\}\right\vert\; \le\;
\displaystyle \frac{1}{k}+\omega_{f,K}\left(\frac{n}{k}\right)
\end{array}
\]
for $i=1,2,\ldots, n$.  The claim  follows by applying the
supremum norm.\epr

Next define
 $\delta(D):=1+c_1+nc_2+c_2\max\{\Vert v\Vert : v\in D\}$.
 It follows  that
 \begin{equation}\label{keyforall}
 G^\eps_{f}[t,x,k/n]\; \; \subseteq\; \;
\{c_1+c_2(||x||+n/k)\}{\cal B}_n\; \; \subseteq\; \;
  \delta(D){\cal B}_n
 \; \; \forall \; t\in [0,T],\; \; x\in D,\; \; k\in \N
.\end{equation}  Theorem \ref{main} now follows from a slight
variant of the argument in \cite{KMW04} which we include as
Appendix \ref{old}.

\subsection{Proof of Theorem \ref{relathm}}
We first make some  general observations that relate our
tangential and Hamiltonian conditions for strong invariance for an
arbitrary multifunction $G:\R^n\rightrightarrows \R^n$. We let
$H_G$ denote the (upper) Hamiltonian for $G$ (see Section
\ref{suff}). We also set $D^o:=\{p\in \R^n: \langle p,d\rangle \le
0\; \forall d\in D\}$  for each  subset $D\subseteq\R^n$; i.e.,
$D^o$ is the {\em polar set} for $D$. Note that $H_G(x,D)\le 0$
holds if and only if $G(x)\subseteq D^o$.
\begin{lemma}
\label{bas} Let $G:\R^n\rightrightarrows \R^n$ and $R$ be a closed
subset of $\R^n$. (i) If $G(x)\subseteq T^B_{\scriptstyle R}(x)$
for all $x\in R$, then $H_G(x, N^P_{\scriptstyle R}(x))\le 0$ for
all $x\in R$. (ii) If $G$ is closed, convex, and nonempty valued
and lower semicontinuous, and if $H_G(x, N^P_{\scriptstyle
R}(x))\le 0$ for all $x\in R$, then $G(x)\subseteq
T^B_{\scriptstyle R}(x)$ for all $x\in R$.
\end{lemma}

\begin{proof}  Since $T^B_{\scriptstyle
R}(x)\subseteq (N^P_{\scriptstyle R}(x))^o$ for all  $x\in R$
(cf. \cite[Exercise II.7.1(d)]{CLSW98} or Appendix A), the
assumptions of (i) imply $G(x)\subseteq (N^P_R(x))^o$ for all
$x\in R$.  This establishes part (i).
 To prove part (ii), first note that its
hypotheses imply $G(x)\subseteq (N^P_{\scriptstyle R}(x))^o$ for
all $x\in R$.  Let $N^C_{\scriptstyle R}$ and  $T^C_{\scriptstyle
R}$ denote the Clarke normal and tangent cones, respectively (cf.
Appendix A for the relevant definitions). We claim that
$G(x)\subseteq (N^C_{\scriptstyle R}(x))^o$ for all $x\in R$. To
verify this claim, fix $x\in R$ and $v\in G(x)$. By Theorem
\ref{mike} in Appendix A, we can find a continuous selection
$s:\R^n\to\R^n$ of $G$ for which $s(x)=v$. Therefore, if $R\ni
x_i\to x$ and $N^P_{\scriptstyle R}(x_i)\ni \zeta_i\to \zeta\in
\R^n$, then $s(x_i)\in G(x_i)$ implies $\langle s(x_i),
\zeta_i\rangle \le 0$ for all $i$. Passing to the limit as $i\to
\infty$  gives $\langle v,\zeta\rangle\le 0$. Therefore, the
characterization (\ref{chara}) for $N^C_R$ in Appendix A  gives
 $v\in
(N^C_{\scriptstyle R}(x))^o$, which proves the claim. It follows
that $G(x)\subseteq (N^C_{\scriptstyle R}(x))^o =T^C_{\scriptstyle
R}(x)\subseteq T^B_{\scriptstyle R}(x)$ for all $x\in R$, which
proves (ii).
\end{proof}

 Returning to the proof of Theorem \ref{relathm}, assume
$F(x,A)\subseteq T^B_S(x)$ for all $x\in S$. We show that $(F,S)$
is strongly invariant  by extending an argument from the appendix
of \cite{AS03} to our more general situation where the dynamics
$F$ may be discontinuous. Let $M$ and $N$ be any compact subsets
of $\R^n$ contained in ${\cal O}$ such that $M\subseteq {\rm
int}(N)$.  Then
 $E:=S\cap N$ is a closed subset of $\R^n$. Let $\phi:\R^n\to
\R$ be any $C^\infty$ function that is identically $1$ on $M$,
strictly positive on ${\rm int}(N)$, and zero elsewhere. Extend
$F$ to $\R^n$ by defining  $F(x,A)\equiv \{0\}$ outside ${\cal O}$
and define the multifunction $F^\sharp:\R^n\times
A\rightrightarrows\R^n$ by $ F^\sharp(x,a)=\phi(x)\, F(x,a)$.  As
observed in  \cite[Lemma A.4, p. 1696]{AS03}, we have:
\begin{lemma}
\label{claimstar} For each $x\in E$, either $F^\sharp(x,A)=\{0\}$
or
 $T^B_E(x)=T^B_S(x)$.
\end{lemma}
 Therefore, since $T^B_E(x)$ is a cone for all $x\in E$, we get
$F^\sharp(x,A)\subseteq T^B_E(x)$ for all $x\in E$. Applying Lemma
\ref{bas}(i)  with $G(x)=F^\sharp(x,A)$  and $R=E$ gives
$H_{F^\sharp}(x,N^P_E(x))\le 0$ for all $x\in E$. Since the
multifunction $x\mapsto F^\sharp(x,A)$ also satisfies
 $(U')$, Corollary \ref{cor1} implies
$(F^\sharp, E)$ is strongly invariant. Therefore, if $T>0$ and
$y:[0,T]\to {\cal O}$ is any trajectory of $F$ such that $y(0)\in
S$, and if we specialize the preceding argument to the compact set
$M=\{y(t): 0\le t\le T\}$, then $y\in {\rm Traj}(F^\sharp)$
(because $\phi\equiv 1$ on $M$) and $y(0)\in E$, so $y$  remains
in $E\subseteq S$. This shows $(F,S)$ is strongly invariant.

Conversely, assume that $(F,S)$ is strongly invariant. We need to
show that $F(x,A)\subseteq T^B_S(x)$ for all $x\in S$. To this
end, fix $\bar x\in S$ and $\bar a\in A$.  Since $S\subseteq {\cal
O}$ and ${\cal O}\subseteq\R^n$ is open, we can find $\mu>0$ such
that $M:=\bar x+\mu {\cal B}_n\subseteq{\cal O}$. Choose any
compact set $N\subseteq {\cal O}$ such that $M\subseteq {\rm
int}(N)$, and let $\phi:\R^n\to \R$ be any $C^\infty$ function
satisfying the requirements above; i.e., $\phi\equiv 1$ on $M$,
$\phi>0$ on ${\rm int}(N)$, and $\phi\equiv 0$ outside  ${\rm
int}(N)$.  Define $F^{\sharp}$ as before. Then $x\mapsto
F^{\sharp}(x,\bar a)$ is closed and convex valued and lower
semicontinuous.

Let $E=S\cap N$ as before, $\zeta\in N^P_E(\bar x)$, and $v\in
F^{\sharp}(\bar x, \bar a)$ be such that
\begin{equation}\label{mix}H_{F^\sharp(\cdot, \bar a)}(\bar x,\zeta)=\langle \zeta,
v\rangle.\end{equation}  Such a $v$ exists because the sets
$D_i(\bar x)$ are compact.
 Reapplying Michael's Selection Theorem (Theorem \ref{mike}
 in Appendix \ref{nonsmo}) provides a continuous selection
 $s:\R^n\to \R^n$ of $x\mapsto F^{\sharp}(x, \bar a)$ such that
$s(\bar x)=v$.  The characterization (\ref{pn}) of the proximal
normal cone in Appendix \ref{nonsmo} gives a constant $\sigma\ge
0$ such that $\langle \zeta, x'-\bar x\rangle \le \sigma ||x'-\bar
x||^2$ for all $x'\in E$. Let $z$ be a $C^1$ local solution of the
initial value problem $\dot z=s(z)$, $z(0)=\bar x$, which we can
assume remains in $M$. Since $\phi\equiv 1$ on $M$, $z$ is also a
trajectory of $F$. Since we are assuming $(F,S)$ is strongly
invariant,  $z$ also stays in $S$. Since $S\cap M\subseteq E$, we
get  $\langle \zeta, z(t)-\bar x\rangle \le \sigma ||z(t)-\bar
x||^2$ for small $t>0$. Since $\dot z(0)=v$, dividing by $t>0$ and
letting $t\to 0^+$ gives
\begin{equation} \label{suping}
\langle \zeta, v\rangle\; \; \le\; \; \sigma \lim_{t\to
0}t||(z(t)-\bar x)/t||^2\; \; =\; \; 0.\end{equation}
 Since $\phi(\bar x)=1$, and since $\zeta\in N^P_E(\bar x)$ and $\bar a\in A$
 were
 arbitrary,
 (\ref{mix})-(\ref{suping})
give $H_{F}(\bar x,N^P_E(\bar x))\le 0$. The preceding argument
applies to any $\bar x\in {\rm int}(N)\cap S$, by choosing $\phi$
to be $1$ near $\bar x$ and compactly supported on $N$, so
$H_{F}(x,N^P_E(x))\le 0$ for all $x\in {\rm int}(N)\cap S$. If we
now fix such a function $\phi$, then we have
$\phi(x)F(x,A)=:F^{\sharp}(x,A)\equiv \{0\}$ for all $x\not\in{\rm
int}(N)$, so $H_{F^{\sharp}}(x,N^P_E(x))\le 0$ for all $x\in N\cap
S=E$.

Applying Lemma \ref{bas}(ii) with $G(x)=F^{\sharp}(x,\bar a)$,
$R=E$, and any $\bar a\in A$ gives $\phi(x)F(x,A)\subseteq
T^B_E(x)$ for all $x\in E$. If $x\in E\cap {\rm int}(N)$ and
$F^{\sharp}(x,A)=\{0\}$, then $F(x,A)=\{0\}$, because $\phi(x)>0$;
while if $x\in E\cap {\rm int}(N)$ and $F^{\sharp}(x,A)\ne \{0\}$,
then Lemma \ref{claimstar}
 and the fact that $\phi(x)>0$ give $F(x,A)\subseteq
T^B_S(x)$.  Therefore, $F(x,A)\subseteq T^B_S(x)$ for all $x\in
E\cap {\rm int}(N)=S\cap {\rm int}(N)$.  Now we fix $\bar x\in S$
and $\mu>0$ such that   $M:=\bar x+\mu {\cal B}_n\subseteq {\cal
O}$ and apply the preceding argument to all possible compact sets
$N\subseteq {\cal O}$ containing $M$ in their interiors. Since
${\cal O}$ is the union of the interiors of such sets $N$, the
preceding argument gives $F(x,A)\subseteq T^B_S(x)$ for all $x\in
S\cap {\cal O}=S$, as desired. This proves Theorem \ref{relathm}.

\section{Infinitesimal Characterizations of Monotonicity}
\label{sec5}

\subsection{Background and Statement of Results}
In this section, we use Theorem \ref{relathm} to prove  new
characterizations of monotonicity for control systems that may be
non-Lipschitz in the state. Monotone control systems were
introduced by Angeli and Sontag in \cite{AS03}  and have since
been applied extensively in systems biology (cf. \cite{AS03,
AS04a, AS04b}). To simplify our exposition, we  only consider
monotone control systems evolving on Euclidean space with input
values in ${\cal U}:={\cal B}_m$ but our results can  be adapted
to systems whose inputs are valued in any ordered Banach space.
The relevant definitions are as follows. We are given two closed
cones $K\subseteq\R^n$ and $K_u\subseteq\R^m$ (called {\em
positivity cones}) which we assume are convex, nonempty, and
pointed (i.e., $K\cap (-K)=\{0\}$ and $K_u\cap (-K_u)=\{0\}$). We
define orders on $\R^n$ and $\R^m$ as follows: $x_1\succeq
x_2\in\R^n$ if and only if $x_1-x_2\in K$, and $u_1\succeq
u_2\in\R^m$ if and only if $u_1-u_2\in K_u$. The ordering on
${\cal U}$ induces an order on the set of controls ${\cal
U}_\infty:=\{{\rm measurable\ } \alpha: [0,\infty)\to {\cal U}\}$
as follows: $\alpha_1\succeq\alpha_2\in {\cal U}_\infty$ if and
only if $\alpha_1(t)-\alpha_2(t)\in K_u$ for a.a. $t\ge 0$.  We
set ${\cal U}^{\scriptstyle [2]}=\{(u_1,u_2)\in {\cal U}\times
{\cal U}: u_1\succeq u_2\}$ and ${\cal U}^{\scriptscriptstyle
[2]}_\infty=\{(\alpha_1,\alpha_2)\in {\cal U}_\infty\times {\cal
U}_\infty: \alpha_1\succeq \alpha_2\}$.

Our main object of study in this section is  the controlled
dynamic
\begin{equation}
\label{mono} \dot x\in  G(x,u):=\prod_{i=1}^ng_i(x,u)D_i(x),\; \;
\; u\in {\cal U},\; \; x\in \tilde X
\end{equation}
evolving on an open set $\tilde X\subseteq\R^n$ in which each
$g_i: \tilde X\times {\cal U}\to\R$ is locally Lipschitz,  $(H_4)$
is satisfied, and the $D_i$'s  take all their values in some
compact set ${\cal D}\subseteq \R$. The dynamic can be viewed as a
locally Lipschitz dynamic  with  non-Lipschitz disturbances $D_i$
acting on its individual components. Since we are assuming $0\in
D_i(x)$ for all $i$ and $x$, the dynamic (\ref{mono}) is clearly
{\em weakly} invariant for any state constraint set $S\subseteq
\tilde X$, since it allows all constant trajectories. However,
since (\ref{mono}) is not necessarily Lipschitz in the state
variable, it may not be {\em strongly} invariant for some state
constraints.
 Following
\cite{AS03}, we also assume there is  a closed set
$X\subseteq\tilde X$ that is the closure of its interior such that
all trajectories of
\begin{equation}
\label{cald} G_{{\cal D}}(x,u)=\prod_{i=1}^ng_i(x,u){\cal D}, \;
\; u\in {\cal U},\; \; x\in \tilde X
\end{equation}
starting in $X$ remain in $X$. We set $\cv={\rm int}(X)$. Since
$G_{{\cal D}}$ is locally Lipschitz in $x$ and ${\cal
D}\subseteq\R$ is compact, the strong invariance of $(G_{{\cal
D}}, X)$ can be checked
 using standard tangent
cone conditions (cf. \cite[Appendix A]{AS03}). On the other hand,
because we allow non-Lipschitz $D_i$'s, (\ref{mono}) is not in
general  tractable by the standard strong invariance
characterizations. We also define $G^{\scriptstyle
[2]}:\R^{2n}\times {\cal U}^{\scriptstyle [2]}\rightrightarrows
\R^{2n}$ by $G^{\scriptstyle [2]}(x,u)=G(x_1,u_1)\times
G(x_2,u_2)$; i.e., two ``stacked''
 copies of $G$ with ordered inputs.

Again following \cite{AS03},  we also consider more general orders
 given by  arbitrary  closed sets $\Gamma\subseteq X\times X$
as follows: We say that $x_1\succeq x_2\in X$ provided
$(x_1,x_2)\in \Gamma$.  This includes the special case of state
spaces ordered by positivity cones $K$ by choosing
$\Gamma=\{(x_1,x_2)\in X \times X: x_1-x_2\in K\}$. With the order
on $\R^n$ expressed  this way, we set $\Gamma_o=\Gamma\cap
(\cv\times \cv)$. We always assume the following approximation
property, which parallels the approximation requirement in
\cite{AS03}:

\begin{itemize}\item[]\begin{itemize}
\item[$(A)$\ \ \ ] For each $\xi\in {\rm bd}(\Gamma)$, $T>0$, and
$y\in {\rm Traj}_T(G^{\scriptstyle [2]},\xi)$, there exist $T'\in
(0,T]$ and two  sequences $\Gamma_o\ni \xi_k\to\xi$ and ${\rm
Traj}_{T'}(G^{\scriptstyle [2]},\xi_k)\ni y_k\to y$ uniformly on
$[0,T']$ as $k\to \infty$.
\end{itemize}\end{itemize}

This agrees with the approximation condition posited in
\cite[p.1686]{AS03} for locally Lipschitz  dynamics $G$,
 since in that case the convergence of $y_k$ to
$y$ follows from continuous dependence on initial values.
 Condition $(A)$ is
satisfied if for example (i) $X$ is convex (or, even more
generally, strictly star-shaped with respect to some interior
point $\xi\in {\cal V}$) and (ii) there exists a neighborhood
${\cal N}$ of ${\rm bd}(X\times X)$ such that $x\mapsto
G^{[2]}(x,{\cal U}^{[2]})$ is Lipschitz on ${\cal N}$. In this
case, the existence of a sequence $\Gamma_o\ni\xi_k\to \xi$ in (A)
for each $\xi\in {\rm bd}(X\times X)$ follows from (i) and the
argument from \cite[p.1686]{AS03}. Given $y\in {\rm
Traj}(G^{[2]},\xi)$, the existence of trajectories ${\rm
Traj}(G^{[2]},\xi_k)\ni y_k\to y$ uniformly in (A)  then follows
because near $\xi$, we can write $G^{[2]}(x,\alpha(t))=J(t,x,{\cal
B}_n)$ for some Lipschitz parametrization $J$ and the input
$\alpha\in {\cal U}^{\scriptscriptstyle [2]}_\infty$ for $y$ (see
\cite[Chapter 9]{AF90}), so $\dot y=J(t,y,\beta)$ for some input
$\beta$ and small times, and then we can apply continuous
dependence on initial conditions to the dynamics $J$. It  suffices
to check $(A)$ for $\xi\in {\rm bd}(X\times X)$, since
$\Gamma\setminus {\rm bd}(X\times X)\subseteq \Gamma_o$.
 On
the other hand, we do not need to assume Lipschitzness of
(\ref{mono}) in a neighborhood of ${\rm bd}(\Gamma_o)$ as was
needed in \cite{AS03}.
 We call
(\ref{mono})  {\em monotone} provided:
\begin{itemize}\item[]\begin{itemize}
\item[$(M)$\ \ \ ]  If $\alpha_1\succeq\alpha_2\in {\cal
U}_\infty$, $x_1\succeq x_2\in X$, $T>0$, $\phi_1\in {\rm
Traj}_T(G,\alpha_1,x_1)$, and $\phi_2\in {\rm
Traj}_T(G,\alpha_2,x_2)$, then $\phi_1(t)\succeq \phi_2(t)$ for
all $t\in [0,T]$.\end{itemize}\end{itemize}

In other words, (\ref{mono}) is monotone provided its flow map
preserves the orders on its inputs and initial states.  Condition
(M) differs slightly from the definition of monotonicity in
\cite{AS03} because our  non-Lipschitz dynamics generally admit
multiple solution trajectories for some choices of inputs and
initial states. Note that  $G$ is monotone if and only if
$(G^{[2]},\Gamma)$ is strongly invariant.  Moreover, $G^{[2]}$
satisfies the hypotheses of Theorem \ref{relathm} (see  Example
\ref{ExhibitA}). For state spaces ordered by positivity cones $K$,
our main result is:

\begin{theorem}
\label{monothm} The dynamics (\ref{mono})  is monotone if and only
if the following condition is satisfied for all $\xi_1,\xi_2\in
{\cal V}$: $(\xi_1\succeq \xi_2$,  $u_1\succeq u_2)$ $\Rightarrow$
$G(\xi_1,u_1)-G(\xi_2,u_2)\subseteq T^B_K(\xi_1-\xi_2)$.
\end{theorem}
It is easy to check that  the tangent cone condition in Theorem
\ref{monothm} is equivalent to the following for each $\xi_1,
\xi_2\in \cv$: $(\xi_1-\xi_2\in {\rm bd}(K), u_1\succeq u_2)$
$\Rightarrow$ $G(\xi_1,u_1)-G(\xi_2,u_2)\subseteq
T^B_K(\xi_1-\xi_2)$. This is because $T^B_K(x)=\R^n$ for all $x\in
{\rm int}(K)$. For orders induced by  closed sets $\Gamma\subseteq
X\times X$, we also prove:
\begin{theorem}
\label{mono2} The dynamics (\ref{mono})  is monotone if and only
if the following condition is satisfied for all $\xi_1,\xi_2\in
{\cal V}$: $(\xi_1\succeq \xi_2$, $u_1\succeq u_2)$ $\Rightarrow$
$G^{\scriptstyle [2]}(\xi,u)\subseteq
T^B_\Gamma(\xi)$.\end{theorem}

 Our
theorems extend the monotonicity characterizations in \cite{AS03}
by  allowing non-Lipschitz dynamics. The dynamics in \cite{AS03}
take the form $G(x,u)=\{g(x,u): u\in {\cal U}\}$ where $g$ is
continuous and locally Lipschitz in $x$ locally uniformly in $u$.
Since we allow each factor $x\mapsto D_i(x)$ in our dynamics
(\ref{mono}) to be set valued and non-Lipschitz, our dynamics may
be non-Lipschitz in the state. Therefore, the strong invariance
theory from \cite[Chapter 4]{CLSW98}  used in \cite{AS03} no
longer applies. Instead, we prove Theorems
\ref{monothm}-\ref{mono2} using our new strong invariance theory
from Section \ref{sec4}.
\subsection{Proof of Theorems \ref{monothm} and \ref{mono2}}

The following lemmas parallel the corresponding lemmas in
\cite[Section III]{AS03}.

\begin{lemma}
\label{L1} $(G,\cv)$ is strongly invariant.
\end{lemma}
\begin{proof} Let $\bar x\in \cv$, $\alpha\in {\cal U}_\infty$, $T>0$, and $y\in
{\rm Traj}_T(G,\alpha,\bar x)$.  We show $y(t)\in \cv$ for all
$t\in [0,T]$. Note that there exist measurable functions
$\beta_i:[0,T]\to\R$ with $\beta_i(t)\in D_i(y(t))$ for all $i$
and a.a. $t\in [0,T]$ such that $y$ is the unique solution of
\begin{equation}
\label{lmd} \dot x=\prod_{i=1}^n g_i(x,\alpha(t))\beta_i(t)
\end{equation} on $[0,T]$ starting at $\bar x$. This follows from
the Filippov Selection Theorem as in Example \ref{ExhibitA}.
 For each $\tilde
x\in \cv$, let $[0,T]\ni t\mapsto z(t,\tilde x)$ be the solution
of (\ref{lmd}) starting at $\tilde x$.  Notice that $z(\cdot,
\tilde x)$ is a trajectory for $G_{{\cal D}}(x,\alpha(t))$ from
(\ref{cald}) for each $\tilde x\in \cv$ and therefore remains in
$X$ by assumption.   Fix $t_o\in [0,T]$. Viewing (\ref{lmd}) as a
system with control inputs $\alpha$ and $\beta_i$ and arguing as
in \cite[Lemma III.6]{AS03} (see \cite[Lemma 4.3.8]{S98}) shows
that the image ${\cal F}({\cal V})$ of the final point map
$\cv\ni\tilde x\mapsto {\cal F}(\tilde x):=z(t_o,\tilde x)\in X$
contains an open neighborhood $W$ of ${\cal F}(\bar x)$.  By
assumption, $W\subseteq {\cal F}({\cal V})\subseteq X$, so
$W\subseteq {\rm int} (X)={\cal V}$.  Therefore, $y(t_o)={\cal
F}(\bar x)\in W\subseteq\cv$. Since $t_o\in [0,T]$ was arbitrary,
$y$ stays in $\cv$.  This proves the lemma.
\end{proof}
Recall that  $G$ is monotone if and only if $(G^{[2]},\Gamma)$ is
strongly invariant.
  We next assume that the order on $\R^n$ is
expressed in terms of  a closed set $\Gamma\subseteq X\times X$ as
above.
\begin{lemma}
\label{lm2a} The dynamics $G$ is monotone if and only if
$(G^{[2]}, \Gamma_o)$ is strongly invariant.
\end{lemma}
\begin{proof}
If $G$ is monotone, then $(G^{[2]}, \Gamma)$ is strongly
invariant. By the previous lemma, $(G^{[2]}, \cv\times \cv)$ is
also strongly invariant, so the sufficiency follows because
$\Gamma_o=\Gamma\cap (\cv\times \cv)$.  Conversely, assume
$(G^{[2]}, \Gamma_o)$ is strongly invariant. We  show $(G^{[2]},
\Gamma)$ is strongly invariant using our approximation hypothesis
$(A)$. Let $T>0$ and $z:[0,T]\to\R^n\times \R^n$ be any trajectory
of $G^{[2]}$ starting at a point in $\Gamma$, and define
\begin{equation}\label{tbardef}\bar t\;
\stackrel{{\rm def}}{=}\; \sup\{t\in [0,T]: z(s)\in \Gamma \,
\forall s\in [0,t]\}.\end{equation} We need to show that $\bar
t=T$.  Suppose the contrary, so $\bar t<T$.  Notice that $z(s)\in
\Gamma$ for all $s\in [0,\bar t]$.  In particular,  $\xi:=z(\bar
t)\in {\rm bd}(\Gamma)$. We apply $(A)$ to the trajectory
$[0,T-\bar t]\ni t\mapsto y(t):=z(t+\bar t)$ of $G^{[2]}$ starting
at $\xi$. This gives $T'\in (0,T-\bar t]$ and sequences
$\Gamma_o\ni \xi_k\to \xi$ and ${\rm Traj}_{T'}(G^{[2]}, \xi_k)\ni
y_k\to y$ uniformly on $[0,T']$. By hypothesis, the $y_k$'s remain
in $\Gamma_o\subseteq \Gamma$, so their uniform limit  also
remains in the closed set $\Gamma$ on $[0,T']$. Therefore, $z$
remains in $\Gamma$ on $[0,\bar t+T']$, which contradicts our
definition of $\bar t$ in (\ref{tbardef}). This establishes that
$(G^{[2]}, \Gamma)$ is strongly invariant and completes the proof
of the lemma.
\end{proof}
We next relate the tangential conditions for  $G$ and $G^{[2]}$
from Theorems \ref{monothm}-\ref{mono2}.
\begin{lemma} \label{related} For any $\xi\in (\xi_1,\xi_2)\in
\Gamma_o$ and $u=(u_1,u_2)\in {\cal U}^{[2]}$, the following three
conditions are equivalent: (a) $G(\xi_1,u_1)-
G(\xi_2,u_2)\subseteq T^B_K(\xi_1-\xi_2)$, (b)
$G^{[2]}(\xi,u)\subseteq T^B_{\Gamma_o}(\xi)$, and (c)
$G^{[2]}(\xi,u)\subseteq T^B_{\Gamma}(\xi)$.
\end{lemma}
\begin{proof}
First fix $\xi=(\xi_1,\xi_2)\in \Gamma_o$, $u=(u_1,u_2)\in {\cal
U}^{[2]}$ and values $v_i\in D_i(\xi_1)$ and $w_i\in D_i(\xi_2)$.
Set
\begin{equation}
\label{newf} f(x,u)=(f_1(x_1,u_1), f_2(x_2,u_2))=\left(
\prod_{i=1}^ng_i(x_1,u_1)v_i, \prod_{i=1}^ng_i(x_2,u_2)w_i\right).
\end{equation}
 Then
$f$ is locally Lipschitz in $x$ although its trajectories  are not
necessarily trajectories of $G^{[2]}$. The lemma  follows if the
following are equivalent for all $\xi\in\Gamma_o$ and $u\in {\cal
U}^{[2]}$:
\begin{itemize}\item[]
($a')$ $f_1(\xi_1,u_1)-f_2(\xi_2,u_2)\in T^B_K(\xi_1-\xi_2)$,
$(b')$ $f(\xi,u)\in T^B_{\Gamma_o}(\xi)$, and $(c')$ $f(\xi,u)\in
T^B_{\Gamma}(\xi)$.\end{itemize} Since (\ref{newf})  is locally
Lipschitz in $x$, this equivalence  follows from the proof of
Lemma III.9 in \cite{AS03}.\end{proof}

We can now prove our two monotonicity theorems.  Assume
(\ref{mono}) is monotone, $(\xi_1,\xi_2)\in \Gamma_o$, and
$u^o=(u^o_1, u^o_2)\in {\cal U}^{[2]}$.  By Lemma \ref{lm2a},
$(G^{[2]}, \Gamma_o)$ is strongly invariant. Applying Theorem
\ref{relathm} from Section \ref{sec4} to the dynamics $\xi\mapsto
G^{[2]}(\xi,u^o)$ in dimension $2n$ with ${\cal O}=\cv\times \cv$
and $S=\Gamma_o$ gives $G^{[2]}(\xi,u^o)\in T^B_{\Gamma_o}(\xi)$
for all $\xi\in \Gamma_o$.  The tangential conditions from both
theorems therefore follow from Lemma \ref{related}. Conversely,
assume either  of these tangential conditions.   By Lemma
\ref{lm2a}, monotonicity of (\ref{mono}) follows once we show that
$(G^{[2]}, \Gamma_o)$ is strongly invariant, but this follows from
Theorem \ref{relathm} applied to $G^{[2]}$ because (b) of Lemma
\ref{related} is satisfied for all $u\in A:={\cal U}^{[2]}$ and
$\xi\in \Gamma_o$. This concludes the proof of our theorems.

\begin{remark}
\label{setdisturb} The tangential conditions from our monotonicity
theorems remain sufficient for monotonicity if we drop the convex
valuedness and lower semicontinuity assumptions on the $D_i$'s and
keep all our other assumptions  the same (by the same proof).
 Moreover, the
requirement that $0\in D_i(x)$ for all $i$ and $x$ can be relaxed
to requiring an increasing sequence $N_1\subseteq
N_2\subseteq\ldots\subseteq N_k\subseteq\ldots$ of compact subsets
of $\cv\times \cv$ such that (i) $\cv\times \cv =\cup_k N_k$ and
(ii) $0\in D_i(x_1)\times D_i(x_2)$ for all $i\in \{1,2,\ldots,
n\}$, all $x=(x_1,x_2)\in {\rm bd}(\Gamma_o\cap N_k)$, and all $k$
(see Section \ref{relclosed}).
\end{remark}

\appendix
\section{Background in Nonsmooth Analysis}
\label{nonsmo}  In this appendix, we  review the necessary
background in nonsmooth analysis; see \cite{CLSW98} for a more
complete treatment and Section \ref{sec2} above for the relevant
notation. Let $S\subseteq \R^n$ be closed and $x\in S$. A vector
$\zeta\in \R^n$ is called a {\em proximal normal}  of $S$ at $x$
provided there exists a constant $\sigma=\sigma(\zeta, x)\ge 0$
such that
\begin{equation}\label{pn}
\langle \zeta,x'-x\rangle \leq \sigma||x'-x||^2.\end{equation} for
all $x'\in S$. The set of all proximal normal vectors of $S$ at
$x$ is denoted by $N^P_S(x)$.  The following {\em local
characterization} of $N^P_S$ also holds: For any $\delta>0$,
$\zeta\in N^P_S(x)$ if and only if there exists
$\sigma=\sigma(\zeta,x)\ge 0$ such that (\ref{pn}) holds for all
$x'\in S\cap (x+\delta {\cal B}_n)$. We also define  the {\em
Clarke tangent cone} $T^C_S$ and the {\em Bouligand ({\rm a.k.a.}
contingent) tangent cone} $T^B_S$ to subsets $S\subseteq\R^n$ as
follows. We say that $v\in T^C_S(x)$ provided for each sequence
$x_i\in S$ converging to $x$ and each sequence $t_i>0$ decreasing
to $0$, there exists a sequence $v_i\to v$ such that $x_i+t_i
v_i\in S$ for all $i$. In particular, if $S=\{0\}$, then
$T^C_S(0)=\{0\}$. The {\em Bouligand tangent cone} to $S$ is
defined by
\[T^B_S(x):=\left\{q\in \R^n: \exists t_i\in (0,\infty) {\rm \ and\ }
 S\ni x_i\to x {\rm\ s.t. \ } t_i\downarrow 0 {\rm \ and\ }
(1/t_i)(x_i-x)\to q\right\}\; \; \forall x\in S.
\]
 Then $T^B_{\scriptstyle S}(x)\subseteq (N^P_{\scriptstyle
S}(x))^o$ for all $x\in S$ and all closed sets $S\subseteq \R^n$,
where the superscript $o$ denotes the polar set; i.e.,
$D^o:=\{p\in \R^n: \langle p, d\rangle \le 0 \, \forall d\in D\}$
for each $D\subseteq\R^n$. The Clarke and Bouligand tangent cones
can differ for some sets $S$ but are known to agree when $S$ is a
closed convex subset of $\R^n$.  We also use the {\em Clarke
normal cone} $N^C_R$ which can be characterized for any closed
subset $R\subseteq\R^n$ as follows:
\begin{equation}\label{chara} \displaystyle N^C_{\scriptstyle
R}(x)=\overline{\rm co}\left\{\lim_{i\to\infty} \zeta_i:
\zeta_i\in N^P_{\scriptstyle R}(x_i), R\ni x_i\to
x\right\}\end{equation} Then $(N^C_R(x))^o=T^C_R(x)\subseteq
T^B_R(x)$ for all $x\in R$.

Next assume $f:\R^n\to (-\infty,\infty]$ is lower semicontinuous
and $x\in{\rm domain} (f):=\{x':f(x')<\infty\}$. Then
$\zeta\in\R^n$ is called a {\em proximal subgradient} for $f$ at
$x$ provided there exist $\sigma>0$ and $\eta>0$ such that
\begin{equation}\label{ps}f(x')\geq f(x) + \langle \zeta,x'-x\rangle - \sigma\Vert
x'-x\Vert ^2\end{equation} for all  $x'\in x+\eta{\cal B}_n$. The
(possibly empty) set of all proximal subgradients for $f$ at $x$
is denoted by $\partial_P f(x)$. When $f$ is the characteristic
function of a closed set $S\subseteq\R^n$ (i.e., $f(x)=0$ if $x\in
S$ and $f(x)=1$ otherwise), it follows from the local
characterization of $N^P_S$ and
 (\ref{ps}) that
$\partial_Pf(x)=N^P_S(x)$ for all $x\in {\rm bd}(S)$ and
$\partial_Pf(x)=\{0\}$ otherwise.

We next state a version of the Clarke-Ledyaev Mean Value
Inequality  (cf. \cite[p. 117]{CLSW98} for its proof). Let $[x,Y]$
denote the closed convex hull of $x\in \R^n$ and $Y\subseteq
\R^n$.

\begin{theorem}\label{MVI}
Assume $x\in \R^n$, $Y\subseteq \R^n$ is compact and convex, and
$\Psi:\R^n\to\R$ is lower semicontinuous.  Then for any
$\delta<\min_{y\in Y}\Psi(y)-\Psi(x)$ and $\lambda>0$, there exist
$z\in [x,Y]+\lambda {\cal B}_n$  and $\zeta\in
\partial_P \Psi(z)$ such that
$\delta<\langle \zeta,y-x\rangle$ for all $y\in Y$.
\end{theorem}

A function $s:\R^n\to\R^n$ is called a {\em selection} of
$F:\R^n\rightrightarrows \R^n$ provided $s(x)\in F(x)$ for all
$x\in \R^n$. The following result is known as Michael's Selection
Theorem (see \cite[Chapter 8]{BL00} or \cite[Corollary
5.59]{RW98}).
\begin{theorem}
\label{mike} Let $F:\R^n\rightrightarrows\R^n$ be lower
semicontinuous and closed, convex, and nonempty valued.  Let $x\in
\R^n$ and $v\in F(x)$.  Then there exists a continuous selection
$s$ of $F$ for which  $s(x)=v$.
\end{theorem}

The following is a variant of the well known \lq\lq  compactness
of trajectories\rq\rq\ lemma which we use in  Appendix \ref{old}.
Its  proof is a special case of the proof of \cite[Theorem
IV.1.11]{CLSW98}.

\begin{lemma}\label{prop2.2.2}Let $\bar
x\in \R^n$, $T>0$, $\tilde f\in \C[0,T]$ be also continuous in
$t$, and $y_{k}:[0,T]\to \R^n$ be a sequence of uniformly bounded
absolutely continuous functions satisfying $y_{k}(0)=\bar x$ for
all $k$. Assume
\begin{equation}
\label{sc}
 \dot  y_{k}(t)\in \tilde f(\tau_k(t), y_{k}(t) + r_k(t)) +
\delta_k(t) {\cal B}_n
\end{equation}
for a.a. $t\in [0,T]$ and all $k$,  where $\delta_k:[0,T]\to
[0,\infty)$ is a sequence of measurable functions that converges
to $0$ in $L^2[0,T]$ as $k\to\infty$, $r_k:[0,T]\to\R$ is a
sequence of measurable functions converging uniformly to $0$ as
$k\to \infty$, and $\tau_k:[0,T]\to[0,\infty)$ is a sequence of
measurable functions converging uniformly to $\tau(t)\equiv t$ as
$k\to \infty$.
  Then there exists a
trajectory $y$ of $\dot y=\tilde f(t,y)$, $y(0)=\bar x$ such that
a subsequence of $y_{k}$ converges to $y$ uniformly on $[0,T]$.
\end{lemma}

We apply Lemma \ref{prop2.2.2} to continuous mollifications of our
feedback realizations $f\in \C[0,T]$ defined as follows. We first
define the {\em standard mollifier}
\[\eta(t)=
\left\{
\begin{array}{ll}
C\, {\rm exp}\left(\frac{1}{t^2-1}\right),& |t|<1\\
0, & |t|\ge 1 \end{array}\right.
\]
where the constant $C>0$ is chosen so that $\int_{\R}\eta(s)ds=1$.
For each $\eps>0$ and $t\in \R$, set
$\eta_\eps(t):=\eta(t/\eps)/\eps$. Notice that
$\int_{\R}\eta_\eps(t)dt=1$ for all $\eps>0$. Define the following
{\em mollifications} of $f\in \C[0,T]$:
\begin{equation}
\label{ftild} f_\eps(t,x):=\int_{\R} f(s,x)\eta_\eps(t-s)ds
\end{equation}
with the convention that $f(s,x)=0$ for all $s\not\in [0,T]$. Then
$f_\eps\in \C[0,T]$ and $(t,x)\mapsto f_\eps(t,x)$ is continuous
for all $\eps>0$. (See \cite[Appendix C]{E98} for the theory of
convolutions and mollifiers.) We  apply Lemma \ref{prop2.2.2} to a
sequence $f_{\eps(k)}$ of mollifications with $\eps(k)>0$
converging to zero using ideas from the usual proof that
$f_{\eps(k)}(\cdot, x)\to f(\cdot, x)$ in $L^1[0,T]$ for each $x$
and $f\in \C[0,T]$ as $k\to \infty$.

\begin{remark}
\label{cr} If $\tau_k(t)\equiv t$ for all $k$ in Lemma
\ref{prop2.2.2}, then the conclusions of the lemma  remain true
even if the $t$-continuity hypothesis on $f\in \C[0,T]$ is
omitted. This follows from the proof of the compactness of
trajectories lemma in \cite{CLSW98} and is  applied to
(\ref{brack}) below.
\end{remark}

\section{Additional Proofs}
\label{old} In this appendix, we summarize the argument from
\cite{KMW04} needed to complete our proof of Theorem \ref{main}.
 We continue to use the notation from Section \ref{4p1} and we set
\begin{equation}
\label{sharp} T':=\min\left\{T,\frac{\gamma}{32\delta(D)}\right\}
\; \; \; {\rm and} \; \; \; h_k:=\frac{\gamma}{32k\delta(D)}
\end{equation}
for all $k\in \N$. Choose $N>2$
 such that
\begin{equation}
\label{fc}
 D+h_k\delta(D){\cal B}_n\; \; \subseteq\; \;  \bar x+\frac{2\gamma}{3}
 {\cal B}_n\; \; \;
 \forall k\ge N.
\end{equation}
  By the choices of $\gamma\in (0,1)$ and
$\delta(D)$,
\begin{equation}
\label{ssc} 0\; <\;   h_k  \; \le\;  \frac{1}{32kg_f[t,x,k/n]}\;
\; \forall t\in [0,T], x\in D, k\in \N.
\end{equation} Next we define $c(k)\equiv {\rm Ceiling}(T'/h_k)$; i.e.,
$c(k)$ is the smallest integer $\ge T'/h_k$.
 For
each $k\ge N$, we then define a partition $\pi(k):  0=
t_{0,k}<t_{1,k}<\cdots< t_{c(k),k}=T'$ by setting
$t_{i,k}=t_{i-1,k}+h_k$ for $i=1,2,\ldots, c(k)-1$.

We next define sequences $x_{0,k},  x_{1,k},  x_{2,k}, \ldots,
x_{c(k),k}$ for $k\ge N$ as follows.  We set $x_{0,k}=\bar x$ and
\begin{equation}\label{ee}x_{1,k}=\bar x+(t_{1,k}-t_{0,k})v_{o,k},
\end{equation} where
$v=v_{o,k}\in G^\eps_{f}[0,\bar x,k/n]$ satisfies the requirement
from Claim \ref{claim1} for the pair $(t,x)=(0,\bar x)$ and
$h=h_k$.
 By (\ref{keyforall}) and (\ref{ee}), we get
\begin{equation}
\label{numb} \Vert x_{1,k}-\bar x\Vert \; \le\;  h_k\delta(D)\; =
\; \frac{\gamma}{32k},
\end{equation}
so $x_{1,k}\in \bar x+\frac{\gamma}{2}{\cal B}_n=D$.  If $c(k)\ge
2$, then we set
  \begin{equation}\label{eee}x_{2,k}=
x_{1,k}+(t_{2,k}-t_{1,k})v_{1,k},\end{equation} where $v_{1,k}\in
G^\eps_{f}[t_{1,k},x_{1,k},k/n]$ satisfies the requirement from
Claim \ref{claim1} for  the pair $(t,x)=(t_{1,k}, x_{1,k})$ and
$h=h_k$. Then (\ref{keyforall}) and (\ref{eee}) give $\Vert
x_{2,k}-x_{1,k}\Vert \le h_k\delta(D)= \frac{\gamma}{32k}$,  so
\begin{equation}\label{sum}
\Vert x_{2,k}-\bar x\Vert \le \frac{\gamma}{16k},\end{equation} by
(\ref{numb}). By (\ref{sum}), $x_{2,k}\in D$.
  We now repeat this process
except with  $x_{2,k}\in D$ replacing $x_{1,k}$. Proceeding
inductively
 gives sequences $v_{i,k}\in G_{f}^\eps[t_{i,k},x_{i,k},k/n]$
 and  $x_{i,k}$ that satisfy
\begin{equation}\label{sts}x_{i+1,k}=x_{i,k}+(t_{i+1,k}-
t_{i,k})v_{i,k}\end{equation} for each index $i=0,1,\ldots,
c(k)-1$. The $v_{i,k}$'s are chosen using Claim \ref{claim1} with
the choices $h=h_k$ and  $(t,x)=(t_{i,k}, x_{i,k})$ for all $i$
and $k$. The choices of $T'$ and $k\ge 2$ and (\ref{sts}) give
\[\Vert x_{i,k}-\bar x\Vert \; \le\;  \frac{c(k)\gamma}{32k} \;
\le\;  \left(\frac{T'}{h_k}+1\right)\frac{\gamma}{4k}\; \le\;
\frac{\gamma}{2}\] for all $i$ and $k$.    It follows that the
sequences $x_{i,k}$ lie in $D$.

For each $k\ge N$,  we then choose $x_{\pi(k)}$ to be the unique
polygonal arc satisfying $x_{\pi(k)}(0)=\bar x$ and
\begin{equation}
\label{poly} \dot
x_{\pi(k)}(t)=f_\eps(\tau_k(t),x_{\pi(k)}(t)+r_k(t))+z_k(\tau_k(t))
\end{equation}
for all $t\in [0,T']\setminus \pi(k)$, where $\tau_k(t)$ is the
partition point  $t_{i,k}\in \pi(k)$ immediately preceding $t$ for
each $t\in [0,T']$,
\begin{equation}
\label{pa} z_k(t_{i,k}):=v_{i,k}-f_\eps(t_{i,k},
x_{\pi(k)}(t_{i,k}))\; \; \forall i,k
\end{equation}
the $v_{i,k}\in G_f^\eps[t_{i,k}, x_{\pi(k)}(t_{i,k}),k/n]$
satisfy the conclusions from  Claim \ref{claim1} for the pairs
$(t,x)=(t_{i,k}, x_{i,k})$ and $h=h_k$, and
\begin{equation}\label{rk}r_k(t):=x_{\pi(k)}(\tau_k(t))-x_{\pi(k)}(t) \; \; \; \forall
t\in [0,T'],\, \forall k.\end{equation} Then $x_{\pi(k)}$ is the
polygonal arc connecting the points $x_{i,k}$ for $i=0,1,2,\ldots,
c(k)$ so $x_{i,k}\equiv x_{\pi(k)}(t_{i,k})$.

 Since $(t,x)\mapsto f_\eps(t,x)$ is continuous, we can
 use Claim \ref{claim2} and the forms of $z_k$ and $r_k$ in
 (\ref{pa})-(\ref{rk}) to
 check  that
 (\ref{poly}) satisfies the requirements of
 our
compactness of trajectories lemma with $y_k=x_{\pi(k)}$ (see
\cite{KMW04}), so we can find a subsequence of $x_{\pi(k)}$ that
converges uniformly to a trajectory $y_\eps$ of $\dot
y=f_\eps(t,y)$, $y(0)=\bar x$. By possibly passing to a
subsequence without relabelling, we can assume that
 $x_{\pi(k)}\to y_\eps$ uniformly on $[0,
 T']$.
 Since $t_{i+1,k}-t_{i,k}\in [0,h_k]$  and
 $x_{i+1,k}=x_{i,k}+(t_{i+1,k}-t_{i,k})v_{i,k}\in D$ for all
$i=0,1,\ldots, c(k)-1$ and $k\ge N$, conditions (\ref{fc}) and
(\ref{ssc}) and  Claim \ref{claim1} give
\begin{equation}\label{kee3}\Psi(x_{i,k})- \Psi(x_{i-1,k})\le
\frac{h_k}{k}\end{equation} for $i=1,2,\ldots, c(k)$. Summing the
inequalities (\ref{kee3}) over $i$ and noting that $h_k\le \gamma$
gives \begin{equation}\label{summm} \Psi(x_{i,k}) \; \; \le\; \;
\Psi(\bar x)+\frac{c(k)h_k}{k}\; \; \le \; \;
 \Psi(\bar x)+\frac{1}{k}(T'+\gamma)\end{equation} for
all $i$ and $k$ . Since $x_{\pi(k)}(\tau_k(t))=x_{i,k}$ for all
$t\in [t_{i,k}, t_{i+1,k})$ and all $i$ and $k$, (\ref{summm})
gives
\begin{equation}
\label{hence} \Psi(x_{\pi(k)}(\tau_k(t)))\le \Psi(\bar
x)+\frac{1}{k}(T'+\gamma)
\end{equation}
for all $t\in [0,T']$.  Since $|\tau_k(t)-t|\le h_k\to 0$ as $k\to
+\infty$ for all  $t\in [0,T']$,  $x_{\pi(k)}(\tau_k(t))\to
y_\eps(t)$ for all $t\in [0,T']$ as $k\to +\infty$. Moreover,
$\Psi$ is lower semicontinuous, so it follows from (\ref{hence})
that
\begin{equation}\label{wee}
\Psi(y_\eps(t))\le \Psi(\bar x)\end{equation} for all $t\in
[0,T']$.

Now let $y_{1/i}:[0,T']\to \R^n$ be the trajectory obtained by the
preceding argument with the choice  $\eps=1/i$ for each $i\in \N$.
Note that $y_{1/i}(t)\in D$ for all $i$ and $t$ because each of
the polygonal arcs $x_{\pi(k)}$ constructed above joins points in
$D$ and $D$ is closed and convex. Moreover,
\begin{equation}
\label{brack}
\begin{array}{lll}
\dot y_{1/i}(t)&=&f(t,y_{1/i}(t))\;  + \;
\left[f_{1/i}(t,y_{1/i}(t))-f(t,y_{1/i}(t))\right]
\end{array}\end{equation} for all $i$ and almost all $t\in [0,T']$.
One can check (see \cite{KMW04}) that the $y_{1/i}$'s are also
equicontinuous, so we can assume (by passing to a subsequence)
there is a continuous function $y:[0,T']\to D$ such that
 $y_{1/i}\to y$ uniformly on   $[0,T']$ (by the
 Ascoli-Arzel\`a lemma).
One can show  (cf. \cite[Claim 4.3]{KMW04}) that
\begin{equation}\label{k3}f_{1/i}(t,y_{1/i}(t))-f(t,y_{1/i}(t))
\to 0 \; \; \; {\rm in}\; \; \;  L^2[0,T'] \end{equation} as $i\to
\infty$.  The proof of (\ref{k3}) is based on the standard proof
that $f_{1/i}(\cdot,x)-f(\cdot,x) \to 0$ in $L^1[0,T']$ for each
$x$  (see for example \cite[pp. 630-1]{E98}) and property $(C_2)$
from the definition of $\C[0,T]$. It therefore follows from Remark
\ref{cr} and the form of the dynamics (\ref{brack}) that a
subsequence of $y_{1/i}$ converges uniformly to some trajectory of
$f$
 on $[0,T']$. This must be the aforementioned function
$y$, which is therefore a trajectory of $f$.
  Again
using the lower semicontinuity of $\Psi$, (\ref{wee}) applied
along the sequence $\eps=1/i$ gives \[\Psi(y(t))\le
\liminf_{i\to\infty}\Psi(y_{1/i}(t))\le \Psi(\bar x)\le 0\] for
all $t\in [0,T']$.  In particular, $y\in {\rm Traj}_{T'}(F,\bar
x)$ and  $y:[0,T']\to {\cal S}$.

Finally, we prove the strong invariance assertion of  the theorem.
We assume $x_o\in {\cal S}$, $T\ge 0$, and $z\in {\rm
Traj}_T(F,x_o)$. Set
\begin{equation}\label{tbardefa}
\bar t := \sup\left\{t\ge 0: \Psi(z(s))\le 0 \, \forall s\in
[0,t]\right\}.\end{equation}  We next show that $\bar t=T$ (by
contradiction), which would imply that
 $z$ remains in ${\cal S}$ on $[0,T]$ and establish the theorem.
 Suppose $\bar t<T$.
 The lower semicontinuity of $\Psi$ gives
$\Psi(z(\bar t\, ))\le 0$. The definition  of $\bar t$  implies
that $z(\cdot)\in {\cal S}$ on $[0,\bar t]$ and in particular
$\bar x:=z(\bar t)\in {\rm bd}({\cal S})\subseteq{\cal U}$. Next
we let $f\in \C'_F([0,T-\bar t], \bar x)$ satisfy the requirement
$(U')$ for $F$ and the following trajectory for $F$:
\begin{equation}\label{ft}[0,T-\bar t]\ni t\mapsto y(t):=z(t+\bar
t).\end{equation} By definition, we can then find
 $\gamma\in (0,1)$
satisfying $f(t,x)\in {\rm ccone}\{F(x)\}$ for a.a. $t\in
[0,\gamma)$ and all $x\in \bar x+\gamma {\cal B}_n$. By reducing
$\gamma>0$ without relabelling, we can assume
 $\bar x+\gamma {\cal B}_n\subseteq {\cal U}$.

 By uniqueness of solutions of the initial value
problem $\dot y=f(t,y)$, $y(0)=z(\bar t)$ on $[0,T-\bar t]$, the
above argument applied to $f$ and $\bar x=z(\bar t)$  gives
$\tilde t\in (0,T-\bar t\, )$ such that
\begin{equation}\label{dip2}
\Psi(z(\bar t+t\, ))\; \; \le\; \; 0 \; \; \; \forall t\in
[0,\tilde t].
\end{equation}
Here we use the fact that the trajectory for $f$ on $[0,T']$
starting at $\bar x$ (and valued in ${\cal S}$) that we
constructed above  (where we can assume $T'<\gamma$) can be
extended to a new trajectory for $f$ defined on all of $[0,T-\bar
t]$ by the linear growth assumption $(C_3)$ on $f$ and so
coincides with (\ref{ft}) by our uniqueness assumption in $(U')$.
We can therefore set $\tilde t=T'$.
 Since $z$ remains in ${\cal S}$ on $[0,\bar t]$, we conclude from
 (\ref{dip2})
 that $z$ remains in ${\cal S}$ on $[0,\bar
t+\tilde t]$. This contradicts the definition  (\ref{tbardefa}) of
$\bar t$ and proves Theorem \ref{main}.

{}\begin{remark} \label{franko} The preceding proof provides a
{\em constructive approach} to finding viable trajectories for our
Carath\'eodory feedback realizations $f$ that remain in ${\cal
S}$.  As suggested by \cite{F05}, an alternative but highly {\em
nonconstructive} proof of Theorem \ref{main} would proceed as
follows.  Let $x(t)$ be any trajectory of $F$ starting in ${\cal
S}$.  By Condition $(U')$, $x(t)$ admits a feedback realization
$f\in {\cal C}[0,T]$, and our Hamiltonian assumption gives
$\langle (f(t,x),0),q\rangle\le 0$ for  all $q\in
N^{\scriptscriptstyle P}_{{\scriptscriptstyle {\rm
epi}(\Psi)}}(x,\Psi(x))=\{(\xi,-1): \xi\in
\partial_P\Psi(x)\}$,
almost all $t\ge 0$, and all $x\in {\cal U}$, where
$\epi(\Psi):=\{(x,r): r\ge \Psi(x)\}$ is the epigraph of $\Psi$.
It is not hard to deduce from this (cf. \cite{CLSW98}) that
$(f(t,x),0)\in T^{\scriptscriptstyle
B}_{{\scriptscriptstyle\epi(\Psi)}}(x,r)$ for almost all $t\ge 0$,
all $x\in {\cal U}$, and all $(x,r)\in {\rm epi}(\Psi)$, where
$T^{\scriptscriptstyle B}$ denotes the Bouligand tangent cone.
Applying the measurable viability theorem (see \cite[Section
4]{FPR95}) to the dynamics $F(t,x)=\{(f(t,x),0)\}$ provides a
trajectory $t\mapsto (\phi(t),\Psi(x(0)))$ for $\dot x=f(t,x)$,
$\dot y=0$ starting at $(x(0),\Psi(x(0)))$ that stays in ${\rm
epi}(\Psi)$. This requires $f$ to be modified outside ${\cal U}$
in the usual way.  Hence, $\Psi(\phi(t))\le\Psi(x(0))\le 0$ for
all $t$.  By the uniqueness part of Condition $(U')$,
$\phi(t)\equiv x(t)$ so $x(t)$ stays in ${\cal S}$. Unfortunately,
the preceding alternative argument is highly nonconstructive,
since the proof of the measurable viability theorem relies on
Zorn's Lemma to construct the trajectory $\phi(t)$.
One natural question which should be considered is how our Euler
constructions from our proof could be used to build numerical
schemes for approximating viable trajectories for Carath\'eodory
dynamics. This type of result could be useful in physical
applications. This question will be addressed by the author in
future research.

\end{remark}


\begin{thebibliography}{99}\addtolength{\itemsep}{-0.5\baselineskip}
\bibitem{AS03}
D. Angeli and E. Sontag, {\em Monotone control systems}. IEEE
Trans. Autom. Control {\bf 48} (2003), 1684--1698.

\bibitem{AS04a}
D. Angeli and E. Sontag, {\em Interconnections of monotone systems
with steady-state characteristics}.  Optimal Control,
Stabilization, and Nonsmooth Analysis (M. de Queiroz, M. Malisoff
and P. Wolenski, Eds.), Lecture Notes in Control and Information
Sciences Vol. 301, Springer-Verlag, Heidelberg, 2004,   135-154.

\bibitem{AS04b}
D. Angeli and E. Sontag, {\em Multistability in monotone
input/output systems.}  Systems and Control Letters {\bf 51}
(2004), 185--202.

\bibitem{AF90}
J.-P. Aubin  and H.  Frankowska, {\em Set-Valued Analysis}.
Systems and Control: Foundations and Applications 2. Birkh\" auser
Boston, Inc., Boston, MA, 1990.

\bibitem{BL00}
J. Borwein and A. Lewis,  {\em Convex Analysis and Nonlinear
Optimization. Theory and Examples}. CMS Books in
Mathematics/Ouvrages de Math\'ematiques de la SMC 3.
Springer-Verlag, New York, 2000.

\bibitem{C75} F. Clarke, {\em Generalized gradients and
applications}. Trans. Amer. Math. Soc. {\bf  205} (1975),
247--262.


\bibitem{CL94}F. Clarke and Y. Ledyaev, {\em Mean value
inequalities in Hilbert space}.  Trans. Amer. Math. Soc. {\bf 344}
(1994), 307--324.



\bibitem{CLR97} F. Clarke, Y.  Ledyaev and M. Radulescu, {\em Approximate
invariance and differential inclusions in Hilbert spaces}. J.
Dynam. Control Syst.  {\bf 3} (1997), 493--518.



\bibitem{CLSW98}
F. Clarke, Y. Ledyaev, R. Stern and P. Wolenski, {\em Nonsmooth
Analysis and Control Theory}. Graduate Texts in Mathematics 178.
Springer-Verlag, New York, 1998.


\bibitem{CS03}
F. Clarke and R. Stern, {\em State constrained feedback
stabilization}.  SIAM J. Control Optim. {\bf 42} (2003), 422--441.


\bibitem{D02} T. Donchev, {\em Properties of the reachable
set of control systems}. Systems and Control Letters {\bf 46}
(2002), 379--386.

\bibitem{DRW04a}
T. Donchev, V. Rios and P. Wolenski, {\em A characterization of
strong invariance for perturbed dissipative systems}.  Optimal
Control, Stabilization, and Nonsmooth Analysis (M. de Queiroz, M.
Malisoff and P. Wolenski, Eds.), Lecture Notes in Control and
Information Sciences Vol. 301, Springer-Verlag, Heidelberg, 2004,
343-349.





\bibitem{E98}
L. Evans, {\em Partial Differential Equations}. Graduate Studies
    in Mathematics, 19. American Mathematical Society, Providence, RI, 1998.
\bibitem{F05}
H. Frankowska, {\em Letter to Michael Malisoff}, March 11, 2005.

\bibitem{FPR95} H. Frankowska, S. Plaskacz and T. Rzezuchowski, {\em Measurable
viability theorems and the Hamilton-Jacobi-Bellman equation.}
Journal of Differential Equations {\bf 116} (1995), 265-305.



\bibitem{KMW04}
M. Krastanov, M. Malisoff and P. Wolenski, {\em On the strong
invariance property for non-Lipschitz dynamics}. LSU Mathematics
Electronic Preprint Series No. 2004-5 and ArXiv math.OC/0403180.


\bibitem{M01}
M. Malisoff, {\em Viscosity solutions of the Bellman equation for
exit time optimal control problems with non--Lipschitz dynamics}.
ESAIM: Control, Optimisation and Calculus of Variations {\bf 6}
(2001), 415--441.


\bibitem{RW03} V. Rios and P. Wolenski, {\em A characterization of
strongly invariant systems for a class of non-Lipschitz
multifunctions}. Proceedings of the  42nd IEEE Conference on
Decision and Control (Maui, HI,  December 2003), 2593-2594.

\bibitem{RW98}
 R. Rockafellar and R. Wets, {\em Variational Analysis.}
 Grundlehren der Mathematischen Wissenschaften
 [Fundamental Principles of Mathematical Sciences] Vol. 317,
 Springer-Verlag, Berlin, 1998.




\bibitem{S98}
E. Sontag, {\em Mathematical Control Theory. Deterministic
Finite-Dimensional Systems. Second Edition}. Texts in Applied
Mathematics 6. Springer-Verlag, New York, 1998.

\bibitem{S03}
H. Sussmann, {\em Uniqueness results for the value function via
direct trajectory-construction methods}. Proceedings of the 42nd
IEEE Conference on Decision and Control (Maui, HI,  December
2003), 3293-3298.

\bibitem{V00}
R. Vinter, {\em Optimal Control}. Systems and  Control:
Foundations and Applications. Birkh\"auser Boston, Inc., Boston,
MA, 2000.

\bibitem{WZ98} P. Wolenski and Y. Zhuang, {\em Proximal analysis and the
minimal time function}.  SIAM J. Control Optim.  {\bf 36} (1998),
1048--1072.

\end{thebibliography}
\end{document}